\documentclass{article}
\usepackage[ansinew]{inputenc}
\usepackage{amsmath}
\usepackage{amssymb}
\usepackage{amsthm}
\usepackage{mathrsfs}
\usepackage[mathscr]{euscript}
\usepackage{latexsym}
\usepackage{cite}
\usepackage{mathpartir}
\usepackage{enumerate}
\usepackage{dsfont}
\usepackage{stmaryrd}
\usepackage{fontenc}
\usepackage{bussproofs}
\usepackage{multicol}
\usepackage{color}
\DeclareMathAlphabet{\mathpzc}{OT1}{pzc}{m}{it}

\newenvironment{calcu}{\begin{array}{rl}}{\end{array}}
\newcommand{\expro}[1]{&#1\\}
\newcommand{\explo}[2]{#1&\ensuremath{\quad\langle\;}\textnormal{#2}\ensuremath{
\;\rangle}\\
}
\newcommand{\cuan}[4]{\ensuremath{(#1 #2 \,|\, #3 \,\boldsymbol{\cdot}\, #4)}}
\newcommand{\cuant}[3]{\ensuremath{(#1 #2 \boldsymbol{\cdot} #3)}}

\theoremstyle{definition}

\makeatletter
\newcommand{\douwidehat}[2]{%
	\sbox0{$\m@th#1\widehat{\hphantom{#2}}\vphantom{t}$}%
	\sbox2{$t$}%
	\dimen2=\ht0
	\advance\dimen2 -\ht2
	\sbox2{$#2$}%
	\dimen0=\ht0
	\rlap{%
		\raisebox{\dimexpr-\dimen0-\dp2-1pt}[0pt][\dimexpr\dimen2+\dp2]{\box0}%
	}
	{#2}%
}
\makeatother

\title{Calculational HoTT}
\author{Bernarda Aldana\footnote{Aldana B. e-mail: bernarda.aldana@escuelaing.edu.co }, Jaime Bohorquez\footnote{Boorquez J. e-mail: jaime.bohorquez@escuelaing.edu.co }, Ernesto Acosta\footnote{Acosta E. e-mail: ernesto.acosta@escuelaing.edu.co }\\
	Escuela Colombiana de Ingenier\'ia\\ Bogot\'a, Colombia}
\date{}

\begin{document}

\maketitle

\begin{abstract}
We found in Homotopy Type Theory (HoTT), a way of representing a first order 
version of intuitionistic logic (ICL, for intuitionistic calculational logic) where, instead of deduction trees, 
corresponding linear calculational formats are used as formal proof-tools; and 
besides this, equality and logical equivalence have preeminence over 
implication. ICL formalisms had been previously adapted by one of the authors 
to intuitionistic logic  from the classical version of the calculational logic 
proposed by Dijkstra and Scholten.
We formally defined \textit{deductive  chains} in HoTT as a representation of the 
linear formats of ICL. Furthermore, we proved using these deductive chains, 
that the equational axioms and rules of ICL have counterparts in HoTT. In doing 
so, we realized that all the induction operators of the basic types in HoTT are 
actually, homotopic equivalences, fact that we proved in this paper. 
Additionally, we propose an informal method to find canonical functions between 
types.
We think that these results could lead to a complete restatement of HoTT where 
equality and homotopic equivalence play a preeminent role. With this approach, 
and by way of calculational methods, effective and elegant formal proofs in HoTT 
are possible through the proposed formal deductive chains by way of appropriate 
formats and  notations.
\end{abstract}

\section{Introduction}
The ability to effectively prove theorems, by both human and mechanical means, is crucial to formal methods. Formal proofs in mathematics and computer science are being studied because they can be verified by a very simple computer program. An open problem in the Computer Mathematics community is the feasibility to fully formalize mathematical proofs \cite{Bar02a}. Here, feasibility is understood as the capability to generate correct formal mathematics with an effort comparable to that of writing a mathematical paper in, say, \LaTeX.

Traditionally, proofs of theorems and formal deductions in deduction systems, are defined, expressed, reasoned about, and performed in principle, through 
formal objects called deduction trees. Typical of these structured forms of 
defining formal deductions are the natural deduction and the sequent systems due 
to Gentzen. Formal deductions are considered too strict and detailed to be used 
in practice by the working mathematician. In fact, except for very short 
proofs, the use of deduction trees gets easily, messy, hard to read and awkward 
to be explained and reasoned about.  

Notwithstanding, for more than thirty years now, a revolution on the way of 
reasoning and proving in mathematics has gained a substantial community of 
enthusiastic practitioners. The \textit{calculational style} of presenting 
proofs introduced by Dijkstra and Scholten \cite{DS90} is a \textit{formal} 
deduction method based on formula manipulation through linear calculational 
formats \cite{vGast90}. This deduction method has been adopted in some books on 
theoretical computer science \cite{Gri93,Bac03,FvG99,Misra01} and appeared in 
papers on set theory, discrete mathematics and combinatorics \cite{AAB17, 
	Boh05, Boh15}.
It was originally devised as an informal but rigorous and practical 
theorem-proving discipline, in which, on one hand, use of equational reasoning 
(understood as mainly based on the preeminence of logical equivalence and 
equalities) is preferred over the traditional one based on logical implication; 
and, on the other hand, the tree-like way of representing formal derivations is 
replaced by what Lifschitz called \textit{calculations} \cite{Lif01}. 
\textit{Calculational logic} and proof methods were formalized for 
\textit{classical} predicate logic by Gries and 
Schneider~\cite{Gri93,DBLP:conf/procomet/gries98} and, subsequently, 
streamlined by Lifschitz \cite{Lif01}. An analogous approach for the case
of intuitionistic predicate logic was developed by one of the authors in 
\cite{Boh08}.

The purpose of this article is to introduce in HoTT a calculational form of reasoning and proving similar to that proposed in  \cite {Boh08} for the intuitionist logic.  In order to formally express HoTT with equality and equivalence playing a preeminent role, we find inspiration in the Curry-Howard isomorphism based on the facts that, on one hand, HoTT is strongly based on the homotopic character of equality and equivalence, and on the other hand, a calculational version of intuitionistic first order logic (ICL) is well established \cite{Boh08}. For this, homotopic equivalence in HoTT plays the role of logical  equivalence in ICL and \textit{deductive chains}, introduced in this work,  play the role of formal \textit{calculations}, term introduced by Lifchitz \cite{Lif01} to formalize Dijkstra and Scholten calculational format. Through this form of reasoning, we could identify judgments in HoTT that represent, under the Curry-Howard isomorphism, the equation rules of the ICL system.  In other words, we want, not only give equivalence a preeminent role in HoTT,  but endow HoTT  with a deduction method based on equational algebraic manipulations that allows for elegant and formal proof constructions, providing a calculational formalization of theorem proving for the case of HoTT by producing (hopefully) human-readable formal proofs based on the linear formats characteristic of the calculational style.

In order to do so, we extend the syntax of type theory introducing an additional judgment that give rise to a conservative extension which facilitates readable proof calculations. We also introduce, as we mention above, an {\it inhabitation format}, that is, a syntactic tool corresponding to the calculational proof format introduced by Dijkstra and Scholten and formalized by Lifchitz with the name of \textit{calculation}  . 

Additionally, we prove the corresponding judgments in HoTT to the basic equational rules in the ICL system. Some of these rules show that induction operators of some of the basic types in HoTT are actually homotopic equivalences, fact that resulted to be true for the rest of induction operators.  

In section 2, we present a brief overview of the main logic principles or rules (algebraic properties, mainly given by equivalences) and notations (Eindhoven quantifiers) used to prove logic theorems calculationally, and the type judgments which correspond, under Curry-Howard isomorphism, to those equational rules. In section 3, we extend HoTT conservatively introducing a new inhabitation judgment which corresponds to a forgetful version of the usual inhabitation judgment, and present some structural rules which will be needed in later sections. In section 4, we define deductive chains as an alternative way of expressing certain derivations of judgments which are sufficient for argumentation in HoTT.
In section 5, we present the basic types of HoTT following the usual four rules:  formation, construction, elimination and computation, but giving the elimination rules a fundamental role as links of deductive chains. In section 6, we introduce the notion of equivalence of types following \cite{hottbook} and study the identification of pairs, functions and natural numbers using deductive chains. Section 7 corresponds to the presentation of the replacement of equivalents by equivalents property of homotopic type-equivalence, which we called Leibniz properties of type-equivalence. In section 8, we prove that all induction operators are actually equivalences, which gives equality and equivalence a preeminent role in HoTT. In section 9, we prove the equational rules stated in section 2 which were not proved in the above sections. In section 10, we present an informal method to find canonical functions between types.

\section{Eindhoven quantifier logic and notation}\label{Eind}

At the \textit{THE} project in Eindhoven, researchers led by E.W. Dijkstra, in 
the 1970's, devised a uniform notation for quantification in first order logic 
and related areas \cite{DS90}. By 
$\cuan{\mathcal{Q}}{x\!:\!T}{range}{term}$\footnote{The original Eindhoven 
style uses colons as separators; the syntax with $|$ and $\boldsymbol{\cdot}$ is one of 
the many subsequent notational variations based on their innovation.} was meant 
that quantifier $\mathcal{Q}$ binds variable $x$ of type $T$ to be constrained 
to satisfy formula \textit{range} within the textual scope delimited by the 
outer parentheses $(...)$, that expression \textit{term} is evaluated for 
each such $x$ and that those values then are combined via an associative and 
commutative operator related to quantifier $\mathcal{Q}$. For brevety, we refer to Eindhoven quantifiers as \textit{operationals}. For the case of 
logical operationals (corresponding to the universal and existential 
quantifiers), the associated operators are respectively, conjunction and 
disjunction considered as binary boolean operations. 

\begin{center}
 $\cuan{\forall}{x\!:\!T}{range}{term}$ means\; for all $x$ in $T$ satisfying 
\textit{range} we have \textit{term},\\
 $\cuan{\exists}{x\!:\!T}{range}{term}$ means\; for some $x$ in $T$ satisfying 
\textit{range} we have \textit{term},
\end{center}
A general shorthand applying to these notations is that an omitted
$|$\textit{range} defaults to $|$\textit{true}. The following so called 
\textit{trade rules} translate these logical notations to the usual first order 
logic formulas\footnote{$\lor$ and $\land$ denote disjunction and conjunction 
respectively, $\Rightarrow$ denote implication and $\equiv$ denotes equivalence. 
If $E$ is a symbolic expression, $E[k/x]$ is the expression obtained by 
replacing every free occurrence of `$x$' in $E$ by `$k$'.}.

\begin{description}
 \item [[Trade\!\!]]\quad $\cuan{\forall}{x\!:\!T}{P}{Q} \equiv 
\cuant{\forall}{x\!:\!T}{P \!\Rightarrow\! Q}$\\
 \hspace*{7mm} $\cuan{\exists}{x\!:\!T}{P}{Q} \equiv \cuant{\exists}{x\!:\!T}{P 
\!\land\! Q}$
\end{description}

The following equational rules (i.e. expressed as logical equivalences) 
correspond to some of the most basic logical axioms and theorems of a 
calculational version of intuitionistic first order logic \cite{Boh08}.

\begin{description}
 \item [[One-Point\!\!]]\quad $\cuan{\forall}{x\!:\!T}{x\!=\!a}{P} 
\equiv P[a/x]$\\
 \hspace*{14mm} $\cuan{\exists}{x\!:\!T}{x\!=\!a}{P} \equiv 
P[a/x]$
 \item [[Equality\!\!]]\quad $\cuan{\forall}{x,y\!:\!T}{x\!=\!y}{P} 
\equiv \cuant{\forall}{x\!:\!T}{P[x/y]}$\\
 \hspace*{11mm} $\cuan{\exists}{x,y\!:\!T}{x\!=\!y}{P} \equiv 
\cuant{\exists}{x\!:\!T}{P[x/y]}$
 \item [[Range Split\!\!]]\quad $\cuan{\forall}{x\!:\!T}{P\lor Q}{R} \equiv 
\cuan{\forall}{x\!:\!T}{P}{R} \land \cuan{\forall}{x\!:\!T}{Q}{R}$\\
 \hspace*{17mm} $\cuan{\exists}{x\!:\!T}{P\lor Q}{R} \equiv 
\cuan{\exists}{x\!:\!T}{P}{R} \lor \cuan{\exists}{x\!:\!T}{Q}{R}$
 \item [[Term Split\!\!]]\quad $\cuan{\forall}{x\!:\!T}{P}{Q\land R} \equiv 
\cuan{\forall}{x\!:\!T}{P}{Q} \land \cuan{\forall}{x\!:\!T}{P}{R}$\\
 \hspace*{15mm} $\cuan{\exists}{x\!:\!T}{P}{Q\lor R} \equiv 
\cuan{\exists}{x\!:\!T}{P}{Q} \lor \cuan{\exists}{x\!:\!T}{P}{R}$
 \item [[Translation\!\!]]\quad $\cuan{\forall}{x\!:\!J}{P}{Q} \equiv 
\cuan{\forall}{y\!:\!K}{P[f(y)/x]}{Q[f(y)/x]}$\\
 \hspace*{16mm} $\cuan{\exists}{x\!:\!J}{P}{Q} \equiv 
\cuan{\exists}{y\!:\!K}{P[f(y)/x]}{Q[f(y)/x]}$\\
where $f$ is a bijection that maps values of type $K$ to values of type $J$.
 \item [[Congruence\!\!]]\quad $\cuan{\forall}{x\!:\!T}{P}{Q\equiv R} 
\Rightarrow (\cuan{\forall}{x\!:\!T}{P}{Q} \equiv 
\cuan{\forall}{x\!:\!T}{P}{R})$\\
 \hspace*{17mm} $\cuan{\forall}{x\!:\!T}{P}{Q\equiv R} \Rightarrow 
(\cuan{\exists}{x\!:\!T}{P}{Q} \equiv \cuan{\exists}{x\!:\!T}{P}{R})$
\item[[Antecedent\!\!]] \quad $R\Rightarrow \cuan{\forall}{x\!:\!T}{P}{Q}\equiv \cuan{\forall}{x\!:\!T}{P}{R\Rightarrow Q}$\\
\hspace*{18mm}$R\Rightarrow \cuan{\exists}{x\!:\!T}{P}{Q}\equiv \cuan{\exists}{x\!:\!T}{P}{R\Rightarrow Q}$\\
when there are not free occurrences of $x$ in $R$. 
\item[[Leibniz principles\!\!]] \quad $ \cuan{\forall}{x,y\!:\!T}{x=y}{f(x)=f(y)}$\\
\hspace*{29.5mm}$\cuan{\exists}{x,y\!:\!T}{x=y}{P(x)\equiv P(y)}$\\
where $f$ is a function that maps values of type $T$ to values of any other type and $P$ is a predicate.
\end{description}
All of these rules have their counterpart in HoTT. In fact, we derive the following judgments which correspond to the above equational rules. In order to write this judgments we have to use the basic types of HoTT and the homotopic equivalence\footnote{The judgment $A\simeq B\!<:$ means that types $A$ and $B$ are quivalent.} that undertakes the role of logical equivalence in logical equational deductions, and the new judgment $A\!<:$ which asserts that $A$ is inhabited without specifying any object.  The definition of homotopic equivalence will be presented in a later section. These are the corresponding rules in HoTT:
\begin{description}
	\item [[One-Point\!\!]]\quad $\prod_{x:A}\prod_{p:x=a}P(x,p) \simeq P(a,\text{refl}_a)\!<:$\\ [0.1cm]
	\hspace*{14mm} $\sum_{x:A}\sum_{p:x=a}P(x,p) \simeq P(a,\text{refl}_a)\!<:$
	\item [[Equality\!\!]]\quad $\prod_{x:A}\prod_{y:A}\prod_{p:x=y}P(x,y,p) \simeq \prod_{x:A}P(x,x,\text{refl}_x)\!<:$\\ [0.1cm]
	\hspace*{11mm} $\sum_{x:A}\sum_{y:A}\sum_{p:x=y}P(x,y,p) \simeq \sum_{x:A}P(x,x,\text{refl}_x)\!<:$
	\item [[Range Split\!\!]]\quad $\prod_{x:A+B} P(x) \simeq \prod_{x:A}P(\text{inl}(x))\times \prod_{x:B}P(\text{inr}(x))\!<:$\\ [0.1cm]
	\hspace*{16mm} $\sum_{x:A+B} P(x) \simeq \sum_{x:A}P(\text{inl}(x))+ \sum_{x:B}P(\text{inr}(x))\!<:$
	\item [[Term Split\!\!]]\quad $\prod_{x:A} (P(x)\times Q(x)) \simeq \prod_{x:A}P(x)\times \prod_{x:A}Q(x)\!<:$\\ [0.1cm]
	\hspace*{16mm} $\sum_{x:A} (P(x)+Q(x)) \simeq \sum_{x:A}P(x)+ \sum_{x:A}Q(x)\!<:$
	\item [[Translation\!\!]]\quad $\prod_{x:A} P(x)\simeq \prod_{y:B}P(g(y))\!<:$\\ [0.1cm]
	\hspace*{16mm} $\sum_{x:A} P(x) \simeq \sum_{y:B} P(g(y))\!<:$\\
	where $g$ is an inhabitant of $B\simeq A$.
	\item [[Congruence\!\!]]\quad $\prod_{x:A} (P(x)\simeq Q(x))\rightarrow  (\prod_{x:A} P(x)\simeq \prod_{x:A}Q(x))\!<:$\\ [0.1cm]
	\hspace*{16mm} $\prod_{x:A} (P(x)\simeq Q(x))\rightarrow  (\sum_{x:A} P(x)\simeq \sum_{x:A}Q(x))\!<:$
	\item[[Antecedent\!\!]] \quad $(R\rightarrow \prod_{x:A}Q(x))\simeq \prod_{x:A}(R\rightarrow Q(x))\!<:$\\ [0.1cm]
	\hspace*{18mm}a) $\sum_{x:A}(R\rightarrow Q(x))\rightarrow (R\rightarrow \sum_{x:A}Q(x))\!<:$\\
	when $R$ does not depend on $x$.\\ [0.1cm]
	\hspace*{18mm}b) $\sum_{x:A}({\mathds 1}\rightarrow Q(x))\simeq ({\mathds 1}\rightarrow \sum_{x:A}Q(x))\!<:$\\
	\item[[Leibniz principles\!\!]] \quad $\prod\limits_{x,y:A}x\!=\!y \rightarrow f(x) \!=\! f(y)\!<:$\\
	\hspace*{29.5mm}$\prod\limits_{x,y:A}x\!=\!y \rightarrow P(x) \!\simeq\! P(y)\!<:$\\
	where $f\!:\!A \rightarrow B$ and $P\!:\!A\to {\cal U}$ is a type family.
\end{description}

A surprising fact about these judgments is that some correspond to homotopic equivalence versions of elimination rules of basic types. In fact, we prove that all elimination rules of the basic types are homotopic equivalences.
  
\section{Extended Syntax of type theory} \label{sec::STT}

In this section we present a formulation of Martin-L\"of theory defining 
terms, judgments and rules of inference inductively in the style of natural 
deduction formalizations. To this formulation, we adjoin an additional 
judgment yielding (by applying its deriving inference rules) a conservative 
extension that allows to perform agile and readable proof calculations.

We suppose the reader is familiar with the syntax of Martin-L\"of type 
theories. and give an overview of the version appearing in \cite{hottbook}.

\subsubsection*{Contexts} \label{subsub::1}

Contexts are finite lists of variable declarations 
$(x_1\!:\!A_1,...,x_n\!:\!A_n)$, for $n\!\geq\!0$, where free variables 
occurring in the $Ai$'s belong to $\{x_1,...,x_{i-1}\}$ when 
$1\!\leq\!i\!\leq\!n$. This list may be empty and indicates that the distinct 
variables $x_1,...,x_n$ are assumed to have types $A_1,...,A_n$, respectively. 
We denote contexts with letters $\Sigma$ and $\Delta$, which may be juxtaposed 
to form larger contexts.

The judgment $\Gamma\; ctx$ formally denotes the fact that $\Gamma$ is a well 
formed context, introduced by the following rules of inference
	\begin{align*}
	& \; \inferrule*[right=ctx-EMP]{ }{\cdot ctx} & &
\inferrule*[right=ctx-EXT]{x_1\!:\!A_1,...,x_{n-1}\!:\!A_{n-1} 
\vdash A_n\!:\!\mathcal{U}_i}{(x_1\!:\!A_1,... , x_n\!:\!A_n)\,ctx}	 
	\end{align*}
with a side condition for the rule \textsc{ctx-EXT}: the variable $x_n$ must be 
distinct from the variables $x_1,...,x_{n-1}$.

\subsubsection*{Forms of judgment} \label{subsub::2}

We first, consider the three usual basic judgments of type theory.
\begin{align*}
  \Gamma\; ctx && \Gamma \vdash a\!:\!A && \Gamma \vdash a\equiv_A a'
\end{align*}
	
$\Gamma\; ctx$ expresses that $\Gamma$ is a (well-formed) context.
$\Gamma \vdash a\!:\!A$ denotes that a term $a$ has (inhabits) type $A$ in 
context $\Gamma$. $\Gamma \vdash a\equiv_A a'$ means that $a$ and $a'$ are 
definitionally equal objects of type $A$ in context $\Gamma$. 

A fourth weaker and derived judgment, the \textit{inhabitation judgment}, will be useful for our purposes:
\[ \Gamma \vdash A\!<: \]
means that the type $A$ is inhabited in context $\Gamma$, that is, for some 
term $a$, judgment $\Gamma \vdash a\!:\!A$ holds. This judgment corresponds to 
a forgetful version of $\Gamma \vdash a\!:\!A$ where the mention of the term 
$a$ inhabiting type $A$ is suppressed.

Since the main inference rule for introducing this judgment is

\[\inferrule*[]{\Gamma \vdash a\!:\!A}{\Gamma \vdash A\!<:}\]
and its remaining derivating inference rules correspond to forgetful versions of derived inference rules from judgments of the form $\Gamma \vdash a\!:\!A$,
this addition only brings forth a conservative extension of the theory.

\subsubsection*{Structural rules} \label{subsub::3}

The following rule expresses that a context holds assumptions, basically by saying that the typing judgments listed in the context may be derived. 

	\[\inferrule*[right=Vble]{(x_1\!:\!A_1,... , x_n\!:\!A_n)\,ctx}{x_1\!:\!A_1,...,x_{n-1}\!:\!A_{n-1} \vdash A_n\!:\!\mathcal{U}_i}\]	 

Although, the following rules corresponding to the principles of \textit{substitution} and \textit{weakening}  are derivable by induction on all possible derivations, we state them. The principles corresponding to typing judgments are given by
	\begin{align*}
	& \; \inferrule*[right=Subst1]{\Gamma \vdash a\!:\!A \\ \Gamma,x\!:\!A,\Delta \vdash b\!:\!B}{\Gamma,\Delta[a/x] \vdash b[a/x]\!:\!B[a/x]} & &
\inferrule*[right=Wkg1]{\Gamma \vdash A\!:\!\mathcal{U}_i \\ \Gamma,\Delta \vdash  b\!:\!B}{\Gamma,x\!:\!A,\Delta \vdash b\!:\!B}	 
	\end{align*}

and the rules for the principles of judgmental (definitional) equality are	
\begin{align*}
& \; \inferrule*[right=Subst2]{\Gamma \vdash a\!:\!A \\ \Gamma,x\!:\!A,\Delta \vdash b\!\equiv_B\!c}{\Gamma,\Delta[a/x] \vdash b[a/x]\!\equiv_{B[a/x]}\!c[a/x]} & &
\inferrule*[right=Wkg2]{\Gamma \vdash A\!:\!\mathcal{U}_i \\ \Gamma,\Delta \vdash  b\equiv_B\!c}{\Gamma,x\!:\!A,\Delta \vdash b\!\equiv_B\!c}	 
\end{align*}	

The following inference rules express the fact that definitional equality is an equivalence relation preserved by typing.
\begin{align*}
& \; \inferrule*[]{\Gamma \vdash a\!:\!A}{\Gamma \vdash a\!\equiv_A\!a} & & \inferrule*[]{\Gamma \vdash a\!\equiv_A\!b}{\Gamma \vdash b\!\equiv_A\!a} & & \inferrule*[right=Tran]{\Gamma \vdash a\!\equiv_A\!b \\ \Gamma \vdash b\!\equiv_A\!c}{\Gamma \vdash a\!\equiv_A\!c} 
\end{align*}
	\begin{align*}
& \; \inferrule*[]{\Gamma \vdash a\!:\!A \\ \Gamma \vdash A\!\equiv\!B\!:\!\mathcal{U}_i}{\Gamma \vdash a\!:\!B} & &
\inferrule*[]{\Gamma \vdash a\!\equiv_A\!b \\ \Gamma \vdash  A\!\equiv\!B\!:\!\mathcal{U}_i}{\Gamma \vdash a\!\equiv_B\!b}	 
	\end{align*}
Besides the inference rule 
\[\inferrule*[right=Inhab]{\Gamma \vdash a\!:\!A}{\Gamma \vdash A\!<:}\]
introducing the inhabitation judgment, we present the following derivating inference rules for this judgment.
	\begin{align*}
& \; \inferrule*[right=Fappl]{\Gamma \vdash A\!<: \\ \Gamma \vdash  A\!\rightarrow\! B\!<:}{\Gamma\vdash B\!<:} & &
\inferrule*[right=Fcomp]{\Gamma \vdash A\!\rightarrow\!B\!<: \\ \Gamma \vdash  B\!\rightarrow\!C\!<:}{\Gamma \vdash A\!\rightarrow\!C\!<:}	 
	\end{align*}
These rules correspond to forgetful versions of the following rules that are easily derived from the original unextended syntax of type theory.
\begin{align*}
& \; \inferrule*[]{\Gamma \vdash a\!:\!A \\ \Gamma \vdash  f\!:\!A\!\rightarrow\! B}{\Gamma\vdash f(a)\!:\!B} & &
\inferrule*[]{\Gamma \vdash f\!:\!A\!\rightarrow\!B \\ \Gamma \vdash  g\!:\!B\!\rightarrow\!C\!}{\Gamma \vdash g\!\circ\!f\!:\!A\!\rightarrow\!C}	 
\end{align*}
An additional structural rule applying definitional equality of types to the 
inhabitation judgment, that we explicitly use, is
\[\inferrule*[right=Tsubs]{\Gamma \vdash A\!<: \\ \Gamma \vdash  
A\!\equiv\! B}{\Gamma\vdash B\!<:}\]

\section{Deductive Chains in Type Theory}

\bigskip
In classical logic, the task is to derive arbitrary valid formulas from a small 
set of axiom schema. In type theory, the basic task is to show that certain type 
can be inhabited from the inhabitation of another types which are related with 
the first through the inference rules introduced before. This will be done by 
means of an {\it inhabitation format}, a syntactic tool that is analogous to the 
calculational proof format introduced by Dijkstra and Scholten \cite{DS90}. \\ [0.1cm]
Before defining an inhabitation format, we present the following  inference rule which can be derived easily from the definition of homotopic equivalence(\cite{hottbook}, (2.4.11), p.79):
\[\inferrule*[right=Heq]{\Gamma \vdash A\simeq B\!<:}{\Gamma\vdash A\to B\!<:},\]
and explicit four of the fairly obvious inference rules, which are used implicitly in type theory most of the time, and correspond to the fact that judgmentally equal things can always be substituted for each other:
\begin{align*}
& \; \inferrule*[right=Repl1l]{\Gamma \vdash A\equiv B}{\Gamma \vdash A\to C\equiv B\to C} & 
& \inferrule*[right=Repl1r]{\Gamma \vdash A\equiv B}{\Gamma \vdash C\to A\equiv C\to B} & 
\end{align*}
\begin{align*}
& \; \inferrule*[right=Repl2l]{\Gamma \vdash A\equiv B}{\Gamma \vdash A\simeq C\equiv B\simeq C} & 
& \inferrule*[right=Repl2r]{\Gamma \vdash A\equiv B}{\Gamma \vdash C\simeq A\equiv C\simeq B} & 
\end{align*}
Given types $A$ and $B$, we temporarily write $A\leadsto B$ to represent the judgments $A\to B\!<:$, the judgment $A\equiv B$ or the judgment $A\simeq B\!<:$. We claim that for all $n\geq 3$, and given a context $\Gamma$, we have the derivation
\[\inferrule*[right=]{\Gamma \vdash A_1\leadsto A_2 \\ \Gamma \vdash A_2\leadsto A_3\\ \cdots\\ \Gamma \vdash A_{n-1}\leadsto A_{n}}{\Gamma \vdash A_1\leadsto A_n}\]
where the conclusion $\Gamma \vdash A_1\leadsto A_n$  corresponds to  
$\Gamma \vdash A_1\to A_n\!<:$ if at least one of the premises is a judgment of the form $\Gamma \vdash A\to B\!<:$, or to $\Gamma \vdash A_1\simeq A_n\!<:$ if none of the premises is of the form  $\Gamma \vdash A\to B\!<:$ and at least one is of the form $\Gamma \vdash A\simeq B\!<:$, or to $\Gamma \vdash A_1\equiv A_n$ if all the premises are of the form $\Gamma \vdash A\equiv B$.\\ [0.1cm] 
We prove our claim by induction. If $n=3$, we have to show that
\[\inferrule*[right=BaseCase]{\Gamma \vdash A_1\leadsto A_2 \\ \Gamma \vdash A_2\leadsto A_3}{\Gamma \vdash A_1\leadsto A_3}\]
Combining the possibilities for $\leadsto$ we have nine cases.\\ [0.1cm]
Cases $(\equiv, \equiv)$, $(\to,\to)$ and $(\simeq,\simeq)$ are \textsc{Tran}, \textsc{Fcomp}, and transitivity of $\simeq$ (\cite{hottbook},Lemma 2.4.12, p. 79), respectively. \\ [0.1cm]
We only derive the first one of the cases $(\to, \equiv)$, $(\equiv, \to)$, $(\simeq, \equiv)$, and  $(\equiv, \simeq)$: 
\[\inferrule*[right=]{\Gamma \vdash A_1\to A_2\!<: \\ \Gamma \vdash A_2\equiv A_3}{\Gamma \vdash A_1\to A_3\!<:},\]
because the rest are derived in the same way. In fact, 
\begin{prooftree}
	\AxiomC{$\Gamma \vdash A_1\to A_2\!<:$}
	\RightLabel{\textsc{Repl1l}}
	\AxiomC{$\Gamma \vdash A_2\equiv A_3$}
	\UnaryInfC{$\Gamma \vdash A_1\to A_2\equiv A_1\to A_3$}
	\RightLabel{\textsc{Tsubs}}
	\BinaryInfC{$\Gamma \vdash A_1\to A_3\!<:$}
\end{prooftree}
From cases $(\to,\simeq)$ and $(\simeq,\to)$ we derive only the first one
\[\inferrule*[right=]{\Gamma \vdash A_1\to A_2\!<: \\ \Gamma \vdash A_2\simeq A_3\!<:}{\Gamma \vdash A_1\to A_3\!<:},\]
the second is done in the same way. In fact,
\begin{prooftree}
	\AxiomC{$\Gamma \vdash A_1\to A_2\!<:$}
	\RightLabel{\textsc{Heq}}
	\AxiomC{$\Gamma \vdash A_2\simeq A_3\!<:$}
	\UnaryInfC{$\Gamma \vdash A_2\to A_3\!<:$}
	\RightLabel{\textsc{Fcomp}}
	\BinaryInfC{$\Gamma \vdash A_1\to A_3\!<:$}
\end{prooftree}
Now, let us suppose that we have the derivation
\[\inferrule*[right=IndHyp]{\Gamma \vdash A_1\leadsto A_2 \\ \Gamma \vdash A_2\leadsto A_3\\ \cdots\\ \Gamma \vdash A_{n-2}\leadsto A_{n-1}}{\Gamma \vdash A_1\leadsto A_{n-1}}.\]
Then,
\begin{prooftree}
	\AxiomC{$\Gamma \vdash A_1\leadsto A_2\cdots \Gamma \vdash A_{n-2}\leadsto A_{n-1}$}
	\RightLabel{\textsc{IndHyp}}
	\UnaryInfC{$\Gamma \vdash A_1\leadsto A_{n-1}$}
	\AxiomC{$\Gamma \vdash A_{n-1}\leadsto A_{n}$}
	\RightLabel{\textsc{BaseCase}}
	\BinaryInfC{$\Gamma \vdash A_1\leadsto A_{n} $}
\end{prooftree}
This proves our claim.\\ [0.1cm]
Due to the rules \textsc{Fappl}, \textsc{Tsubs} and \textsc{Heq}  we have the derivation
\[\inferrule*[right=]{\Gamma \vdash a:A \\ \Gamma \vdash A\leadsto B}{\Gamma \vdash B\!<:}.\] 
Let us suppose a given context $\Gamma$. A \textit{deductive chain} is a derivation of the form 
\begin{equation}\label{Deriva}
\inferrule*[right=]{\overset{\vdots}{\Gamma \vdash a:A_1} \\ \overset{\vdots}{\Gamma \vdash A_1\leadsto A_2} \\ \cdots\\ \overset{\vdots}{\Gamma \vdash A_{n-1}\leadsto A_{n}}}{\Gamma \vdash A_{n}\!<:}.
\end{equation} represented schematically as a vertical deductive chain:
\[
\begin{array}{rl}
 & A_n \\ \leftrightarrows & \\ & A_{n-1} \\ & \vdots \\  \leftrightarrows &\\ & A_2 \\   \leftrightarrows & \\ & A_1 \\ \stackrel{\mbox{\tiny $\wedge$}}{\mbox{\tiny :}} & \left\langle \textit{inhabitation statement} \right\rangle\\ & a
\end{array}
\]
 These chains, and their concrete versions, will be referred as \textit{inhabitation formats}.
Each link 
\[
\begin{array}{rl}
& B \\ \leftrightarrows & \\ & A 
\end{array}
\]
in the above format, corresponds to one of the  following concrete versions:
\[
\begin{array}{l}
\phantom{\leftarrow } B \\ \leftarrow \left\langle\!:\,\, ; \textit{statement of inhabitation} \right\rangle \\ \phantom{\leftarrow } A 
\end{array}
\]
called \textit{consequence link}, 
\[
\begin{array}{l}
\phantom{\leftarrow } B \\ \equiv \left\langle \textit{evidence of equivalence} \right\rangle  \\ \phantom{\leftarrow } A 
\end{array}\]
called \textit{equivalence link}, or 
\[
\begin{array}{l}
\phantom{\leftarrow } B \\ \simeq \left\langle\!:\,\, ; \textit{statement of inhabitation} \right\rangle   \\ \phantom{\leftarrow } A 
\end{array}
\]
called \textit{homotopic equivalence link}. The closing link, that is the link at th bottom of the deduction chain, 
\[
\begin{array}{rl}
 & A\\ \stackrel{\mbox{\tiny $\wedge$}}{\mbox{\tiny :}} & \left\langle \textit{inhabitation statement} \right\rangle\\ & a
\end{array}
\]
is called \textit{inhabitation link}.\\ [0.1cm] 
In short, this inhabitation format is a deductive chain that represents the concatenation of the premises of a derivation of the form (\ref{Deriva}). Each link of the chain is a judgment of the form $A\!\rightarrow\! B\!<:$, $A\!\equiv\!B\!$, $A\!\simeq\!B\!$ or $a\!:\!A$ written vertically, together with an evidence or a statement supporting it, which is written between angular parentheses. \\ [0.1cm]
If $f:A\rightarrow B$, $g:B\rightarrow C$, 
$h:A\rightarrow B$ and $a\!:\!A$ then $h(g(f(a)))\!:\!D$. This detailed account 
of inhabitation is represented by the following chain:
\[
\begin{calcu}
\expro{D}
\\
\explo{\leftarrow}{\!\!:\,\,$h$}
\\
\expro{C}
\\
\explo{ \leftarrow}{\!\!:\,\,$g$}
\\
\expro{B}
\\
\explo{ \leftarrow}{\!\!:\,\,$f$}
\\
\expro{A}
\explo{\stackrel{\mbox{\tiny $\wedge$}}{\mbox{\tiny :}}}{evidence of inhabitation}
\\
\expro{a.}
\end{calcu}
\]
that derives, not only that $D$ is inhabited, but that $D$ is inhabited by $h(g(f(a)))$. \\ [0.1cm]
Before illustrating the use of deduction chains we introduce some basic 
types in order to present some consequence links which come with their 
specifications.

\section{Basic Types}\label{DistArr}
We follow the general pattern for introducing new types in Type Theory 
presented in the  HoTT book \cite{hottbook}. The specification of a type 
consist mainly in four steps: (i)\textit{Formation rules}, (ii) \textit{Construction rules},
(iii) \textit{Elimination rules}, and \textit{Computation rules}. Here, we express the elimination rules in terms of consequence links.\\ [0.1cm]
We assign a special Greek letter to each induction operator introduced in the respective elimination rule. Namely
\[\begin{array}{|c|c|c|c|c|c|c|} \hline
\text{Type} & \Sigma & + & \mathbb N & = & \mathds O & \mathds 1\\ \hline
\text{Induction operator}& \boldsymbol{\sigma} & \boldsymbol{\kappa} & \boldsymbol{\nu} & \boldsymbol{\iota} & \boldsymbol{o} & \boldsymbol{\mu} \\ \hline
\end{array}
\]
\textbf{$\Pi$-types}. The dependent function types or $\Pi$-types, are the most fundamental basic 
types and its elimination rule does not provide links for deductive chains.   \\ 
[0.3 cm]
Given types $A\!:\!{\cal U}$ and $B\!:\!A\rightarrow{\cal U}$ we form the type 
$\prod_{x:A}B(x)\!:\!\mathcal{U}$. For $b\!:\!B$ we construct 
$\lambda(x\!:\!A).b$ of type 
$\prod_{x:A}B(x)$.\\ [0.1cm]
For $f\!:\!\prod_{x:A}B(x)$ and $a\!:\!A$ then $f(a)\!:\!B[a/x]$ and the 
computation rule is
\[(\lambda(x:A).b)(a)\equiv b[a/x]\]
When $B$ does not depend on the objects of $A$, the product type is the function 
type $A\rightarrow B$:
\[
\prod\limits_{x:A}B(x)\; \equiv \; A\rightarrow B.
\]
The propositional reading of $f:\prod_{x: A} B (x)$ is that $f$  is a proof 
that all objects of type $A$ satisfy the property $ B $. We  use this {\it 
	semantic} throughout the paper as necessary. By the way, the elimination 
rules 
of $\Sigma$-types, co-product types, $\mathbb N$-type, and $W$-types, establish
that to prove that all objects of these types satisfy a property,
you have to prove that their constructed objects satisfy the property, and for 
this, the rule introduces an induction operator fulfilling that task.\\ [0.1cm]
One useful property of $\Pi$ types is \textsl{$\Pi$-distribution over arrows}. Let us suppose that for each $x\!:\!A$ we have a function $\varphi_x:P(x)\to 
Q(x)$. Then we can define the  function
$$\Delta : (\prod_{x:A}P(x))\to 
(\prod_{x:A}Q(x))$$ by
$
\Delta(u)(x):\equiv  \varphi_x(u(x)).
$
This shows that if $\prod_{x:A}(P(x)\to Q(x))<:$ then $(\prod_{x:A}P(x))\to 
\prod_{x:A}Q(x)<:$. This property is known as $\Pi$-distribution over arrows 
and 
is frequently used in deductive chains as the following consequence link
\begin{equation}\label{DistArrow}
\begin{calcu}
\expro{\prod_{x:A}Q(x)}
\\
\explo{\leftarrow}{\!:\,$\Delta$\,;\, Definition of $\varphi_x$}
\\
\expro{\prod_{x:A}P(x)}
\end{calcu}
\end{equation}
Later, in the section \ref{InhArr}, we explain a method to find 
definitions of functions such as the one for $\Delta$.\\ [0.1cm]

\textbf{$\Sigma$-types}.The dependent pair types or $\Sigma$-types, are the types whose inhabitants are 
dependent pairs.\\ [0.3 cm]
Given $A\!:\!\mathcal{U}$ and $B\!:\!A\rightarrow\mathcal{U}$ we form 
$\sum_{x:A}B(x)\!:\!\mathcal{U}$ and if $a\!:\!A$ and $b\!:\!B[x/a]$ then 
$(a,b)\!:\! \sum_{x:A}B(x)$.\\ [0.1cm] 
In order to prove a property $C:\sum_{x:A}B(x)\rightarrow 
{\cal U}$ for all objects of the $\Sigma$-type, i.e., to inhabit 
$\prod_{p:\sum_{x:A}B(x)}C(p)$, we must prove the property for its constructed 
objects, 
i.e., to inhabit $\prod_{x:A}\prod_{y:B(x)}C((x,y))$ For this there is a 
function $\boldsymbol{\sigma}(C)$ carrying a proof $g$ of this latter expression 
to the proof $\boldsymbol{\sigma}(C)(g)$ of the former expression. Therefore,
the elimination rule is given by the following consequence link
\[
\begin{calcu}
\expro{\prod\limits_{p:\sum_{x:A}B(x)}C(p)}
\\
\explo{\leftarrow}{\!:\,\,$\boldsymbol{\sigma}_C$}
\\
\expro{\prod\limits_{x:A}\prod\limits_{y:B(x)}C((x,y))}
\end{calcu}
\]
The computation rule states the definition of the function 
$\boldsymbol{\sigma}_C$:
\[
\boldsymbol{\sigma}_C(g)((a,b))\equiv g(a)(b).
\]
For the case when $C$  is a constant family, we have that the induction operator link reduces to 
\[
\begin{calcu}
\expro{(\sum\limits_{x:A}B(x))\to C}
\\
\explo{\leftarrow}{\!:\,\,$\boldsymbol{\sigma}_C$}
\\
\expro{\prod\limits_{x:A}(B(x)\to C)}
\end{calcu}
\]
With the induction operator we can also define functions on $\Sigma$-types. For 
instance, projection functions $\text{pr}_1$ and $\text{pr}_2$ are defined by 
\[\text{pr}_1 :\equiv \boldsymbol{\sigma}_A (g)\;\;\text{and\;\;}\text{pr}_2 :\equiv \boldsymbol{\sigma}_{B\circ \text{pr}_1}(h),\]
where $g:\equiv \lambda(x:A).\lambda(y:B(x)).x$, and  $h:\equiv \lambda(x:A).\lambda(y:B(x)).y$.\\ [0.1cm]
When $B$ does not depend on the objects of $A$, the $\Sigma$-type is the type 
$A\times B$, the Cartesian product type of $A$ and $B$:
\[
\sum_{x:A}B(x)\; \equiv \; A\times B.
\]
\textbf{Coproduct types}. The coproduct corresponds to the disjoint union of sets in Set Theory. \\ [0.1cm]
Given $A\!:\!\mathcal{U}$ and $B\!:\!\mathcal{U}$ we form $A+B\!:\!\mathcal{U}$ 
and if $a\!:\!A$ and $b\!:\!B$ then $\text{inl}(a)\!:A+B$ and 
$\text{inr}(b)\!:\!A+B$.\\ [0.1cm] 
In order to prove a property $C:A+B\rightarrow 
{\cal U}$ for all objects of the coproduct type, i.e., to inhabit 
$\prod_{p:A+B}C(p)$, we must prove the property for its constructed objects, 
i.e., to inhabit $\prod_{x:A}C(\text{inl}(x))\times\prod_{y:B}C(\text{inr}(y))$. 
For this there is a 
function $\boldsymbol{\kappa}_C$ carrying a proof $g$ of the latter type to the 
proof $\boldsymbol{\kappa}_C(g)$ of the former one. Therefore, the elimination 
rule is given by the following consequence link
\[
\begin{calcu}
\expro{\prod\limits_{p:A+B}C(p)}
\\
\explo{\leftarrow}{\!:\,\,$\boldsymbol{\kappa}_C$}
\\
\expro{\prod\limits_{x:A}C(\text{inl}(x))\times\prod\limits_{y:B}C(\text{inr}
(y))}
\end{calcu}
\]
The computation rule states the definition of the function $\boldsymbol{\kappa}_C$:
\[\boldsymbol{\kappa}_C(g)(\text{inl}(a)):\equiv (\text{pr}_1g)(a)\;\; \text{and}\;\; \boldsymbol{\kappa}_C(g)(\text{inr}(b)):\equiv (\text{pr}_2g)(b)\]
\textbf{Empty type}. It is presented as $\mathds O$. This type has no objects 
and its elimination rule is given by the function 
\[
\boldsymbol{o}_C : \prod\limits_{x:{\mathds O}}C(x),
\] 
which states that all the objects of ${\mathds O}$ satisfy any property $C:\mathds O\to {\cal U}$, and there is no computation rule.\\ [0.1cm]
\textbf{Unit type}. It is presented as ${\mathds 1}$. This type has just one 
object, its constructor is $*\!:\!\mathds 1$, and  its elimination rule is 
given by the following link:
\[
\begin{calcu}
\expro{\prod\limits_{x:\mathds 1}C(x)}
\\
\explo{\leftarrow}{\!:\, $\boldsymbol{\mu}_C$}
\\
\expro{C(*)}
\end{calcu}
\] 
which states that in order to prove a property $C:\mathds 1\to {\cal U}$ it is enough to inhabit $C(*)$. Its computation rule is $\boldsymbol{\mu}_C(u)(x):\equiv u$.\\ [0.1cm]
\textbf{The type of natural numbers} is presented as $\mathbb N$ and its constructors are $0\!:\!\mathbb N$ and $s\!:\!\mathbb 
N\rightarrow\mathbb N$.\\ [0.1cm] 
In order to prove a property $C:\mathbb N\rightarrow 
{\cal U}$ for all objects of $\mathbb N$, i.e., to inhabit 
$\prod_{p:\mathbb N}C(p)$, we must prove the property for its constructed 
objects, 
i.e., to inhabit $C(0)\times \left( \prod_{p:\mathbb N}C(p)\rightarrow 
C(s(p))\right) $. For this, there is a 
function $\boldsymbol{\nu}_C$ carrying a proof $g$ of the latter type to the 
proof $\boldsymbol{\nu}_C(g)$ of the former one. Therefore, the elimination rule 
is given by the following consequence link
\[
\begin{calcu}
\expro{\prod\limits_{p:\mathbb N}C(p)}
\\
\explo{\leftarrow}{\!:\,\,$\boldsymbol{\nu}_C$}
\\
\expro{C(0)\times \prod\limits_{p:\mathbb N}C(p)\rightarrow C(s(p))}
\end{calcu}
\]
The computation rule states the definition of the function 
$\boldsymbol{\nu}_C$:
\[
\boldsymbol{\nu}_C(g)(0)\equiv (\text{pr}_1g)(0)
\;\; \text{and}\;\; \boldsymbol{\nu}_C(g)(s(p)))\equiv (\text{pr}_2g)(p,\boldsymbol{\nu}_C(g)(p)).
\]

\textbf{Identity type}. Given any pair of objects $a$ and $b$ of a type $P:{\cal U}$, there is a type 
$(a=_{\phantom{.}_P}\!b):{\cal U}$, called identity type. There is only one 
constructor:
\[
\text{refl}: \prod\limits_{x:P} (x=_{\phantom{.}_P}\!x) 
\]
that states de identification of an object with itself. The objects of $x=y$ are 
called paths from $x$ to $y$.\\ [0.3 cm]
In order to prove a property $C:\prod_{x,y:P}x=y\rightarrow 
{\cal U}$ for all objects of the identity type, i.e., to inhabit 
$\prod_{x,y}\prod_{p:x=y}C(p)$, we must prove the property for its constructed 
objects, 
i.e., to inhabit $\prod_{x:P}C(\text{refl}_x)$. For this there is a 
function $\boldsymbol{\iota}_C$ carrying a proof $g$ of the latter type to the 
proof $\boldsymbol{\iota}_C(g)$ of the former one. Therefore,
the elimination rule is given by the following consequence link
\[
\begin{calcu}
\expro{\prod\limits_{x,y:P}\prod\limits_{p:x=y}C(x,y,p)}
\\
\explo{\leftarrow}{\!:\,\,$\boldsymbol{\iota}_C$}
\\
\expro{\prod\limits_{x:P}C(x,x,\text{refl}_x)}
\end{calcu}
\]
The computation rule states the definition of the function 
$\boldsymbol{\iota}_C$:
\[
\boldsymbol{\iota}_C(g)(x,x, \text{refl}_x):\equiv g(x).
\]
\textbf{Remark}. Induction operators depend on a type family; however, the corresponding computation rules do not. Recall that computation rules for $\boldsymbol{\sigma}$, $\boldsymbol{\kappa}$, $\boldsymbol{\iota}$ and $\boldsymbol{\mu}$, for example, are respectively:
	$\boldsymbol{\sigma}(u)((x,y)):\equiv u(x)(y)$, $\boldsymbol{\kappa}(u,v)(\text{inl}(x)):\equiv u(x)$, $\boldsymbol{\kappa}(u,v)(\text{inr}(y)):\equiv v(y)$,
	$\boldsymbol{\iota}(u)(x,x,\text{refl}_x):\equiv u(x)$, and
	$\boldsymbol{\mu}(u)(*):\equiv u$.
	These computations are idependent of the family type to which they apply. From now on, we do not mention the type families to which they apply .\\ [0.1cm]
	
With the identity induction operator, one can characterize the inhabitants of 
Cartesian product types and coproduct types, this allows us to present the first 
examples of deductive chains.
For the case of the Cartesian product type, 
if $A$ and $B$ are types, then  
\begin{equation}\label{UniqPairs}
\prod\limits_{u:A\times B}u=(\text{pr}_1(u),\text{pr}_2(u))<:
\end{equation}
In fact, 
\[
\begin{calcu}
\expro{\prod\limits_{u:A\times B}u=(\text{pr}_1(u),\text{pr}_2(u))}
\\
\explo{\leftarrow}{\!:\,$\boldsymbol{\sigma}$}
\\
\expro{\prod\limits_{x:A}\prod\limits_{y:B} (x,y)=(\text{pr}_1((x,y)),\text{pr}_2((x,y)))}
\\
\explo{\equiv}{Definition of $\text{pr}_1$ and $\text{pr}_2$}
\\
\expro{\prod\limits_{x:A}\prod\limits_{y:B} (x,y)=(x,y)}
\\
\explo{\stackrel{\mbox{\tiny $\wedge$}}{\mbox{\tiny :}}}{$h(x)(y):\equiv \text{refl}_{(x,y)}$}
\\
\expro{h.}
\end{calcu}
\]
And, for the case of the coproduct type, if $A$ and $B$ are types, then 
 \[\prod\limits_{p:A+B}(\sum\limits_{x:A}(p=\text{inl}(x))+\sum\limits_{y:B}(p=\text{inr}(y)))<:\]
In fact,
\[
\begin{calcu}
\expro{\prod\limits_{p:A+B}\sum\limits_{x:A}(p=\text{inl}(x))+\sum\limits_{y:B}p=\text{inr}(y)}
\\
\explo{\leftarrow}{\!:$\boldsymbol{\kappa}$}
\\
\expro{\phantom{\times}\prod\limits_{a:A}(\sum\limits_{x:A}(\text{inl}(a)=\text{inl}
(x))+\sum\limits_{y:B}\text{inl}(a)=\text{inr}(y))}
\\
&\times\prod\limits_{b:B}\sum\limits_{x:A}(\text{inl}(b)=\text{inl}(x))+\sum\limits_{y:B}
\text {inr}(b)=\text{inr}(y)
\\
\\
\explo{\leftarrow}{\!:\,$\varphi\,\,;\,\, \varphi(u,v):\equiv (\text{inl}\circ 
u, \text{inr}\circ v )$ }
\\
\expro{\prod\limits_{a:A}(\sum\limits_{x:A}\text{inl}(a)=\text{inl}(x)) 
\,\times\, \prod\limits_{b:B}\sum\limits_{y:B}\text{inr}(b)=\text{inr}(y)}
\\
\explo{\stackrel{\mbox{\tiny $\wedge$}}{\mbox{\tiny :}}}{$h:\equiv 
(\lambda a.(a,\text{refl}_{\text{inl}(a)}),\lambda 
b.(b,\text{refl}_{\text{inr}(b)}))$}
\\
\expro{h}
\end{calcu}
\]
\section{Equivalence of types}
Now, we introduce the notion of equivalence of types, but first, we need the 
one 
of 
homotopic functions. Details of this topic may be found in  \cite{hottbook}. \\ 
[0.3 cm]
Let $f$ and $g$ be two dependent functions inhabiting  $\prod_{x:A}P(x)$. We 
say that $f$ and $g$ are homotopic if the type $f\sim g$ defined by
\[
f \sim g :\equiv \prod\limits_{x:A} (f(x)=g(x)) 
\]
is inhabited. 
Two types $A$ and $B$ are equivalent if there is a function $f:A\rightarrow B$ 
such that the type isequiv($f$) defined by 
\[
\text{isequiv}(f) :\equiv ( \sum\limits_{g:B \to A} f\circ g \sim 
\text{id}_B) \times ( \sum\limits_{h:B \to A} h\circ f \sim 
\text{id}_A )  
\]
is inhabited. Therefore, $A$ and $B$ are equivalent if the type $A\simeq B$ 
defined by $\sum_{f:A\to B}\text{isequiv}(f)$ is inhabited. However, in order 
to 
prove equivalence in this paper, we do not use the type isequiv($f$), but the 
type qinv($f$), which is a simpler equivalent  version (see \cite{hottbook}, 
2.4 
p. 76) and is defined by  
\[
\text{qinv}(f) :\equiv  \sum\limits_{g:B \to A} \left( (f\circ g \sim 
\text{id}_B) \times (g\circ f \sim \text{id}_A) \right) . 
\]
This means that in order to show that types $A$ and $B$ are equivalent we must 
exhibit a 4-tuple 
\[
\boldsymbol{f}:\equiv (f,f',\alpha, \alpha')
\]
where 
\[ f:A\to B, \quad f':B\to A,\quad \alpha: f\circ f' \sim \text{id}_B, \quad  
\text{and}\quad  \alpha': f'\circ f \sim \text{id}_A.\]
For instance, let us show that given types $A$ and $B$, 
\begin{equation}\label{comm+}
A+B\simeq B+A<:
\end{equation}
 In fact, let  $f\!:\!A+B\to B+A$ and $f'\!:\!B+A\to A+B$ be defined by\linebreak 
$f(\text{inl}(a)):\equiv  \text{inr}(a)$, $f(\text{inr}(b)):\equiv  
\text{inl}(b)$,  $f'(\text{inl}(b)):\equiv  \text{inr}(b)$ and 
$f'(\text{inr}(a)):\equiv  \text{inl}(a)$. Then, the folowing deductive chain 
shows that  $f\circ f'\sim \text{id}_{B+A}$ is inhabited:  
\[
\begin{calcu}
\expro{f\circ f'\sim \text{id}_{B+A}}
\\
\explo{\equiv}{Definition of $\sim$}
\\
\expro{\prod\limits_{p:B+A} f(f'(p))= p}
\\
\explo{\leftarrow}{\!:\,$\boldsymbol{\kappa}$}
\\
\expro{\prod\limits_{b:B} (f(f'(\text{inl}(b)))= \text{inl}(b)) \times 
\prod\limits_{a:A} (f(f'(\text{inr}(a)))= \text{inr}(a))}
\\
\explo{\equiv }{Definition of $f$ and $f'$}
\\
\expro{\prod\limits_{b:B} (\text{inl}(b)= \text{inl}(b)) \times 
\prod\limits_{a:A} (\text{inr}(a)= \text{inr}(a))}
\\
\explo{\stackrel{\mbox{\tiny $\wedge$}}{\mbox{\tiny :}} }{$u:\equiv \lambda b.\text{refl}_{\text{inl}(b)}\,\,;\,\, 
v:\equiv \lambda a.\text{refl}_{\text{inr}(a)}$}
\\
\expro{(u,v)}
\end{calcu}
\]
We prove $f'\circ f \sim \text{id}_{A+B}<:$ in the same way.\\ [0.1cm]
We present three equivalences characterizing the identification of objects of 
certain types: pairs, functions, and natural numbers. \\ 
[0.1cm]
\textbf{Identification of pairs}. Let $A$, $B$ be types. Then for all $u$ and $v$ inhabitants of $A\times B$ we 
have that 
\[
u=v\;\simeq\; (\text{pr}_1(u)=\text{pr}_1(v))\times 
(\text{pr}_2(u)=\text{pr}_2(v))\, <: 
\]
\textit{Proof}. First of all, we define $P_1(u,v):\equiv \text{pr}_1(u)=\text{pr}_1(v)$ and $P_2(u,v):\equiv \text{pr}_2(u)=\text{pr}_2(v)$. And now, we define $f\!:\!u\!=\!v\to P_1(u,v)\times P_2(u,v)$, by means of the following deductive chain:
\[
\begin{calcu}
\expro{\prod\limits_{u,v:A\times B} \prod\limits_{p:u=v}P_1(u,v)\times P_2(u,v)}
\\
\explo{\leftarrow}{\!:\,$\boldsymbol{\iota_1}$}
\\
\expro{\prod\limits_{u:A\times B} P_1(u,u)\times P_2(u,u)}
\\
\explo{\stackrel{\mbox{\tiny $\wedge$}}{\mbox{\tiny :}} }{$h:\equiv \lambda u.(\text{refl}_{\text{pr}_1(u)},\text{refl}_{\text{pr}_2(u)}$}
\\
\expro{h}
\end{calcu}
\] 
Therefore we may define $f:\equiv \boldsymbol{\iota_1}(h)(u,v)$.\\ [0.1cm]
In order to define a function $f':P_1(u,v)\times P_2(u,v) \to  u\!=\!v$, let us consider the following 
deductive chain:
\[
\begin{calcu}
\expro{\prod\limits_{u,v:A\times B}P_1(u,v)\times P_2(u,v)\to u\!=\!v}
\\
\explo{\leftarrow}{\!:\,$\boldsymbol{\sigma}\,\,;\,\,\boldsymbol{\sigma}(w)((a,
	c),(b,d),(p,q)):\equiv w(a)(b)(c)(d)(p)(q)$}
\\
\expro{\prod\limits_{a,b:A}\prod\limits_{c,d:B}\prod\limits_{p:a=b}\prod\limits_
	{q:c=d}(a,c)=(b,d)}
\\
\explo{\leftarrow}{\!:\,$\boldsymbol{\iota_2}\,\,;\,\,\boldsymbol{\iota_2}(z)(a,
a , c ,
	c,\text{refl}_a,\text{refl}_c):\equiv z(a)(c)$}
\\
\expro{\prod\limits_{a:A}\prod\limits_{c:B}(a,c)=(a,c)}

\\
\explo{\stackrel{\mbox{\tiny $\wedge$}}{\mbox{\tiny :}} }{$k(a,c):\equiv \text{refl}_{(a,c)}$}
\\
\expro{k}
\end{calcu}
\] 
Therefore, we can put $f':\equiv 
(\boldsymbol{\sigma}\!\circ\!\boldsymbol{\iota_2})(k)(u,v)$.\\ [0.3 cm]
Now, let us show that $\prod_{u,v:A\times B}f\circ f'\sim \text{id}<:$
\small{\[
\begin{calcu}
\expro{\prod\limits_{u,v:A\times 
		B}\prod\limits_{g:P_1(u,v)\times P_2(u,v)}f(f'(g))=g}
\\
\explo{\equiv}{Definition of $f$ and $f'$}
\\
\expro{\prod\limits_{u,v:A\times 
		B}\prod\limits_{g:P_1(u,v)\times P_2(u,v)}
		
(\boldsymbol{\iota_1}(h)(u,v))\left((\boldsymbol{\sigma}
\!\circ\!\boldsymbol{\iota_2})(k)(u,v))(p,q)\right)=(p,q)}
\\
\explo{\leftarrow}{\!:\,$\boldsymbol{\sigma}$}
\\
\expro{\prod\limits_{a,b:A}\prod\limits_{c,d:B}\prod\limits_{p:a=b}\prod\limits_
	{q:c=d}
(\boldsymbol{\iota_1}(h)((a,c),(b,d)))\left((\boldsymbol{\sigma}
\!\circ\!\boldsymbol{\iota_2})(k)((a,c),(b,d)) (p,q)\right)=(p,q)}
\\
\explo{\leftarrow}{\!:$\boldsymbol{\iota}$}
\\
\expro{\prod\limits_{a:A}\prod\limits_{c:B}
(\boldsymbol{\iota_1}(h)((a,c),(a,c)))\left((\boldsymbol{\sigma}
\!\circ\!\boldsymbol{\iota_2})(k)((a,c),(a,c)) 
(\text{refl}_a,\text{refl}_c)\right)=(\text{refl}_a,\text{refl}_c)}
\\
\explo{\equiv}{Definition of $\boldsymbol{\sigma}$, $\boldsymbol{\iota_2}$, and 
	$k$}
\\
\expro{\prod\limits_{a:A}\prod\limits_{c:B}(\boldsymbol{\iota_1}(h)((a,c),(a,
	c)))(\text{refl}_{(a,c)})=(\text{refl}_a,\text{refl}_c)}
\\
\explo{\equiv}{Definition of $\boldsymbol{\iota_1}$, and $h$}
\\
\expro{\prod\limits_{a:A}\prod\limits_{c:B}(\text{refl}_a,\text{refl}_c)=(\text{refl}_a,\text{refl}_c)}
\\
\explo{\stackrel{\mbox{\tiny $\wedge$}}{\mbox{\tiny :}} }{$j:\equiv \lambda a.\lambda c.\text{refl}_{(\text{refl}_a,\text{refl}_c)}$}
\\
\expro{j}
\end{calcu}
\] }
The proof of $\prod_{u,v:A\times B}f'\circ f\sim \text{id}<:$ is done in the same way.\\ [0.1cm]
As a particular case, we have that if $a, c: A$, and $b, d: B$, then
\begin{equation}
(a,b)=(c,d)\; \simeq \; a\!=\!c\times b\!=\!d <:
\end{equation}
\textbf{Identification of functions}. Let $A$ and $B$ be two types, and $f$ and $g$ objects of $A\to B$. Then 
\begin{equation}\label{FuncExt}
f=g\; \simeq \; f\sim g \,<:
\end{equation}
The inhabitation can not be proved with the theory introduced till now but 
introduced as an axiom in \cite{hottbook} as {\it function 
extensionality}.\\[2mm]
\textbf{Identification of natural numbers}. If one introduces the type family 
\[
\text{code}: \mathbb N \to \mathbb N \to {\cal U}
\]
defined by
\[\text{code}(0,0):\equiv  \mathds 1, \; \text{code}(s(n),0):\equiv  \mathds 
O,\;\text{code}(0,s(n)):\equiv \mathds O,\; \text{and}\]
\[\text{code}(s(m),s(n)):\equiv  \text{code}(m,n)\]
then, theorem 2.13.1 in \cite{hottbook} states that, for all $m,n:\mathbb N$, 
we have that 
\begin{equation}
m=n\; \simeq\; \text{code}(m,n)<:
\end{equation}
Its proof introduces the functions $\textit{encode}\!\!:\prod_{m,n:\mathbb N} 
m\!=\!n\to \text{code}(m,n)$ and $\textit{decode}\!:\!\prod_{m,n:\mathbb N} 
\text{code}(m,n)\to m=n$, and shows that the functions $\text{encode}(m,n)$ and 
$\text{decode}(m,n)$ are q-inverses of each other.\\ [0.1cm]
In next sections, we explore several properties related with equivalence.

\section{Leibniz properties of type equivalence}

By Leibniz properties, we refer to the replacement of equivalents by 
equivalents (or congruence) property of, in this case, homotopic 
type-equivalence. 
\subsection{Leibniz principles.}

These are precisely [\textbf{Leibniz principles}] mentioned in section 2, and refer to the fact that equality  is preserved respectively, by function application and type dependency (through, equivalence)\\ [0.1cm]
Let $A,B:{\cal U}$, $f\!:\!A \rightarrow B$ and $P\!:\!A\to {\cal U}$. Then
 \[\prod\limits_{x,y:A}x\!=\!y \rightarrow f(x) \!=\! f(y)<:\quad\text{and}\quad
\prod\limits_{x,y:A}x\!=\!y \rightarrow P(x) \!\simeq\! P(y)<:\]
In fact,
\[
\begin{calcu}
\expro{\prod\limits_{x,y:A}\prod\limits_{p:x\!=\!y} f(x) \!=\! f(y)}
\\
\explo{\simeq}{\!:\! $\boldsymbol{\iota}$}
\\
\expro{\prod\limits_{x:A} f(x) \!=\! f(x)}
\\
\explo{\stackrel{\mbox{\tiny $\wedge$}}{\mbox{\tiny :}}}{$h(x):\equiv \text{refl}_{f(x)}$}
\\
\expro{h}
\end{calcu}
\]
One defines $\text{ap}_{f}(x,y,p):\equiv \boldsymbol{\iota}(h)(x,y,p)$, and by 
definition of $\boldsymbol{\iota}$, we get\\ 
$\text{ap}_{f}(x,x,\text{refl}_x):\equiv \boldsymbol{\iota} 
(h)(x,x,\text{refl}_x):\equiv h(x):\equiv \text{refl}_{f(x)}$.\\ [0.1cm]On the 
other hand,
\[
\begin{calcu}
\expro{\prod\limits_{x,y:A}\prod\limits_{p:x\!=\!y} P(x) \!\simeq\! P(y)}
\\
\explo{\simeq}{\!:\! $\boldsymbol{\iota}$}
\\
\expro{\prod\limits_{x:A} P(x) \!\simeq\! P(x)}
\\
\explo{\stackrel{\mbox{\tiny $\wedge$}}{\mbox{\tiny :}}}{$k(x):\equiv \text{id}_{P(x)}$}
\\
\expro{k}
\end{calcu}
\]
One defines $\text{tr}^{P}(x,y,p)\!\!:\equiv \!\!\boldsymbol{\iota}(k)(x,y,p)$\footnote{This object is called transport$^P$ in the HoTT book \cite{hottbook} }, and by 
definition of $\boldsymbol{\iota}$, we get\linebreak
$\text{tr}^{P}(x,x,\text{refl}_x):\equiv 
\boldsymbol{\iota}(k)(x,x,\text{refl}_x):\equiv k(x):\equiv \text{id}_{P(x)}$.

\subsection{Leibniz inference rules.}\label{LeibInf}

Leibniz inference rules generally express the fact that  type equivalence is preserved by replacement, in any given type expression, of any of its subexpressions by an equivalent one. We derive Leibniz inference rules for coproduct types, and for $\Pi$ and $\Sigma$ types, which are precisely [\textbf{Congruence}] and [\textbf{Translation}] rules, endowing HoTT, by this means, with a calculational style of proof.\\ [0.1cm]
Let $A,B.C:{\cal U}$ and $P,Q: A\to {\cal U}$. Then
\begin{description}
	\item[[Congruence\!\!]]\quad
	$\prod_{x:A} P(x)\simeq Q(x)\to\prod_{x:A} P(x)\simeq \prod_{x:A} Q(x)\!<:$ \quad $\Pi${\sc Eq}1\\ [0.3cm]
	\hspace*{16mm} $\prod_{x:A} P(x)\simeq Q(x)\to\sum_{x:A} P(x)\simeq \sum_{x:A} Q(x)\!<:$\quad $\Sigma${\sc Eq}1
\vspace{0.3cm}
	\item [[Translation\!\!]]\quad 
	$\prod_{\boldsymbol{f}:A\simeq B}\left( \prod_{x:A} P(x)\simeq \prod_{y:B} P(f'(y))\right) <:$\quad $\Pi${\sc Eq}2 \\ [0.3cm] 
	\hspace*{16mm} $\prod_{\boldsymbol{f}:A\simeq B}\left( \sum_{x:A} P(x)\simeq \sum_{y:B} P(f'(y))\right) <:$\quad $\Sigma${\sc Eq}2
\vspace{0.3cm}
	\item[[Coproduct Monotony]\!\!] \quad $(A\simeq B)\to (A+C\simeq B+C)<:$\quad $+${\sc Eq}1\\ [0.3cm]
	\hspace*{34mm}$(A\simeq B) \to (C+A\simeq C+B)<:$ \quad $+${\sc Eq}2
\end{description}  
\vspace{0.3cm}
\textit{Proof of $\Pi${\sc Eq}1.}
Suppose that $\boldsymbol{\Phi}\!:\!\prod_{x:A}P(x)\simeq Q(x)$, with 
$\boldsymbol{\Phi}(x)\equiv (\phi_x,\phi_x', \alpha, \alpha')$,  
$\alpha \!:\!\phi_x\!\circ\! \phi_x' \sim \text{id}_{Q(x)}$ and  $\alpha': 
\phi_x'\!\circ\! \phi_x \sim \text{id}_{P(x)}$.  Let 
\[
\psi:\prod_{x:A}P(x)\to \prod_{x:A}Q(x)
\]
be defined by $\psi(f)(x):\equiv \phi_x(f(x))$\footnote{$\psi$ is precisely the function $\Delta$ of $\Pi$-distibution over arrows, see (\ref{DistArrow}) } and let 
\[\psi':\prod_{x:A}Q(x)\to \prod_{x:A}P(x)
\]
be defined by $\psi'(g)(x):\equiv  \phi_x' (g(x))$. Observe that
\begin{equation}\label{CalPeq1}
\psi(\psi'(g))(x)
\equiv
\phi_x(\psi'(g)(x))
\equiv
\phi_x(\phi_x'(g(x)))
\equiv
(\phi_x\circ \phi_x')(g(x))
\end{equation}

Then, in order to prove $\psi\circ \psi'\sim \text{id}<:$\,, it is enough to prove $(\psi\circ\psi')(g)=g<:$\,. for all $g:\prod_{x:A}Q(x)$. In fact,
\[
\begin{calcu}
\expro{(\psi\circ\psi')(g)=g}
\\
\explo{\simeq}{Function extensionality (\ref{FuncExt})}
\\
\expro{\prod\limits_{x:A}(\psi\circ\psi')(g)(x)= 
g(x)}
\\
\explo{\equiv}{See above calculations (\ref{CalPeq1})}
\\
\expro{\prod\limits_{x:A}(\phi_x\circ 
\phi_x')(g(x))= g(x)}
\\
\explo{\stackrel{\mbox{\tiny $\wedge$}}{\mbox{\tiny :}}}{$u(g)(x):\equiv \alpha (g(x))$}
\\
\expro{u}
\end{calcu}
\]
The proof of $\psi\circ \psi'\sim \text{id}<:$ is done similarly.\\ [0.1cm]
\textit{Proof of $\Sigma${\sc Eq}1.}
Suppose that $\boldsymbol{\Phi}:\prod_{x:A}P(x)\simeq Q(x)$ with 
$\boldsymbol{\Phi}(x)\equiv (\phi_x,\phi_x', \alpha, \alpha')$, $\alpha 
:\phi_x\circ \phi_x' \sim \text{id}_{Q(x)}$, and  $\alpha': \phi_x'\circ \phi_x 
\sim \text{id}_{P(x)}$.  Let 
\[
\psi:\sum_{x:A}P(x)\rightarrow \sum_{x:A}Q(x),
\] 
be defined by  $\psi(p):\equiv 
(\text{pr}_1(p),\phi_{\text{pr}_1(p)}(\text{pr}_2(p))$ and let 
\[
\psi':\sum_{x:A}Q(x)\rightarrow \sum_{x:A}P(x)
\]
be defined by $\psi'(q):\equiv 
(\text{pr}_1(q),\phi_{\text{pr}_1(q)}'(\text{pr}_2(q))) $. Observe that 
\begin{equation}\label{CalSeq1}
\psi(\psi'((x,y)))
\equiv
\psi((x,\phi_x'(y)))
\equiv
(x,\phi_x(\phi_x'(y)))
\equiv
(x,(\phi_x\circ \phi_x')(y))
\end{equation}
Then,
\[
\begin{calcu}
\expro{\psi\circ \psi'\sim \text{id}}
\\
\explo{\equiv}{Definition of $\sim$}
\\
\expro{\prod\limits_{q:\sum_{x:A}Q(x)}(\psi\circ\psi')(q)=q}
\\
\explo{\leftarrow}{\!:\,$\boldsymbol{\sigma}$}
\\
\expro{\prod\limits_{x:A}\prod\limits_{y:Q(x)}\psi(\psi'((x,y)))=(x,y)}
\\
\explo{\equiv}{See above computations (\ref{CalSeq1})}
\\
\expro{\prod\limits_{x:A}\prod\limits_{y:Q(x)}(x,(\phi\circ\phi')(y))=(x,y)}
\\
\explo{\simeq}{$(a,b)=(c,d)\simeq a=c\times b=d\, <:$  \,;\, $\Pi${\sc Eq}1}
\\
\expro{\prod\limits_{x:A}\prod\limits_{y:Q(x)}x=x\times (\phi\circ\phi')(y)=y}
\\
\explo{\stackrel{\mbox{\tiny $\wedge$}}{\mbox{\tiny :}}}{$h(x,y):\equiv 
(\text{refl}_x,\alpha(y))\,\,;\,\,\alpha:\phi\circ \phi'\sim \text{id}$}
\\
\expro{h}
\end{calcu}
\]
We prove $\psi'\circ \psi\sim \text{id}<:$ similarly.\\ [0.1cm]
\textit{Proof of $\Pi${\sc eq2.}}  
Suppose that $\boldsymbol{f}:A\simeq B$. Let  
\[\psi:\prod_{x:A} P(x)\rightarrow \prod_{y:B}P(f'(y))
\]
be defined by $\psi(u)(y):\equiv u(f'(y))$, and let 
\[
\psi':\prod_{y:B} P(f'(y))\rightarrow \prod_{x:A}P(x)
\]
be defined by $\psi'(v)(x):\equiv v(f(x))$.  Let us see that  $\psi'$ is a
quasi-inverse of $\psi$. On one hand, we have
\[
\begin{calcu}
\expro{\psi\circ \psi'\sim \text{id}}
\\
\explo{\equiv}{Definition of $\sim $}
\\
\expro{\prod\limits_{v:\prod_{y:B} P(f'(y))}\psi(\psi'(v))=v}
\\
\explo{\equiv}{Definition of $\psi$ and $\psi'$}
\\
\expro{\prod\limits_{v:\prod_{y:B} P(f'(y))}v\circ f \circ f'=v}
\\
\explo{\simeq}{Function extensionality (\ref{FuncExt})\,\,;\, $\Pi${\sc Eq}1}
\\
\expro{\prod\limits _{v:\prod_{y:B} P(f'(y))}v\circ f \circ f'\sim v}
\\
\explo{\equiv}{Definition of $\sim$}
\\
\expro{\prod\limits_{v:\prod_{y:B} P(f'(y))}\prod\limits_{y:B}v(f( f'(y)))=v(y)}
\\
\explo{\leftarrow}{\!:\,$\Delta$ \,;\,  
$\varphi_{(v,y)}:\equiv \text{ap}_v(f( f'(y)),y)$,\, see (\ref{DistArrow}) }
\\
\expro{\prod\limits_{v:\prod_{y:B} P(f'(y))}\prod\limits_{y:B}f( f'(y))=y}
\\
\explo{\leftarrow}{\!:\, $\lambda z.(\lambda v.z)$}
\\
\expro{\prod\limits_{y:B}f(f'(y))=y}
\\
\explo{\equiv}{Definition of $\sim$}
\\
\expro{f\circ f'\sim \text{id}_B}
\\
\explo{\stackrel{\mbox{\tiny $\wedge$}}{\mbox{\tiny :}}}{Hypothesis}
\\
\expro{\alpha}
\end{calcu}
\]
On the other hand, we can show, exactly in the same way, that 
\[h'\circ h\sim \text{id}_{\prod_{x:A} P(x)}<:.\]
\textbf{Application of $\Pi$-translation rule} (to prove  
$\text{isSet}(\mathbb N)<:$).\\[1mm]
We can use the translation rule to prove  $\text{isSet}(\mathbb 
N)<:$\footnote{See definition 3.1.1 in \cite{hottbook}}\,. In 
fact, let  $\Phi:m=n\rightarrow \text{code}(m,n)$ be defined by $\Phi:\equiv 
\text{encode}(m,n)$ and let $\Psi:\text{code}(m,n)\rightarrow m=n$ be defined by 
$\Psi:\equiv \text{decode}(m,n)$. Then,
\[
\begin{calcu}
\expro{\text{isSet}(\mathbb N)}
\\
\explo{\equiv}{Definition of isSet}
\\
\expro{\prod\limits_{m,n:\mathbb N}\prod\limits_{p,q:m=n}p=q}
\\
\explo{\simeq}{$\Pi$-translation rule\,;\,$m=n\simeq\text{code}(m,n)$}
\\
\expro{\prod\limits_{m,n:\mathbb 
		N}\prod\limits_{s,t:\text{code}(m,n)}\Psi(s)=\Psi(t)}
\\
\explo{\stackrel{\mbox{\tiny $\wedge$}}{\mbox{\tiny :}}}{See definition of $h$ below}
\\
\expro{h}
\end{calcu}
\]
where $h$ is defined by
\[
h(m,n,s,t)=
\begin{cases}
\boldsymbol{\mu}_{1}(\boldsymbol{\mu}_{2}(\text{refl}_{\Psi(*)}))&\text{ if } 
\text{code}(m,n)={\mathds 
	1}\\
\boldsymbol{o}_{C}(s)(t), &\text{ if } \text{code}(m,n)={\mathds O}
\end{cases}
\]
with $C\equiv \prod\limits_{t:\textbf{0}}\Psi(s)=\Psi(t)$. 
The definition of $h$ is justified by
\[
\begin{calcu}
\expro{\prod\limits_{s,t:\textbf{1}}\Psi(s)=\Psi(t)}
\\
\explo{\leftarrow}{\!:\,$\boldsymbol{\mu}_{1}$} 
\\
\expro{\prod\limits_{t:{\mathds1}}\Psi(\ast)=\Psi(t)}
\\
\explo{\leftarrow}{\!:\,$\boldsymbol{\mu}_{2}$}
\\
\expro{\Psi(\ast)=\Psi(\ast)}
\\
\explo{\stackrel{\mbox{\tiny $\wedge$}}{\mbox{\tiny :}}}{$u:\equiv \text{refl}_{\Psi(\ast)}$}
\\
\expro{u}
\end{calcu}
\]
\textit{Proof of $\Sigma${\sc eq2.}}
Suppose that $\boldsymbol f:A\simeq B$. Let 
\[
\psi:\sum_{x:A}P(x)\rightarrow \sum_{y:B}P(f'(y))
\]
defined by  $\psi(u):\equiv (f(\text{pr}_1(u)),\text{pr}_2(u))$ and let 
\[
\psi':\sum_{y:B}P(f'(y))\rightarrow \sum_{x:A}P(x)
\]
defined by $\psi'(v):\equiv (f'(\text{pr}_1(v)),\text{pr}_2(v))$. 
Observe that 
\begin{equation}\label{CalSeq2} 
\psi(\psi(v))\equiv\psi((f'(\text{pr}_1(v)),\text{pr}_2(v))\equiv((f\circ 
f')(\text{pr}_1(v)),\text{pr}_2(v))\end{equation}
Then we have that
\[
\begin{calcu}
\expro{\psi\circ \psi'\sim \text{id}}
\\
\explo{\equiv}{Definition of $\sim$}
\\
\expro{\prod\limits_{v:\sum_{y:B}P(f'(y))}\psi(\psi'(v))=v}
\\
\explo{\leftarrow}{\!:\,$\boldsymbol{\sigma}$}
\\
\expro{\prod\limits_{y:B}\prod\limits_{z:P(f'(y))}\psi(\psi'(y,z))=(y,z)}
\end{calcu}
\]
\[
\begin{calcu}
\explo{\equiv}{See above calculations (\ref{CalSeq2})}
\\
\expro{\prod\limits_{y:B}\prod\limits_{z:P(f'(y))}((f\circ f')(y),z)=(y,z)}
\\
\explo{\simeq}{$(a,b)=(c,d)\simeq (a\!=\!c)\times (b\!=\!d)\,<:$ \,;\, 
$\Pi${\sc Eq}1}
\\
\expro{\prod\limits_{y:B}\prod\limits_{z:P(f'(y))}((f\circ f')(y)\!=\!y)\times 
(z\!=\!z)}
\\
\explo{\stackrel{\mbox{\tiny $\wedge$}}{\mbox{\tiny :}}}{$h(y,z):\equiv (\alpha(y),\text{refl}_z)$}
\\
\expro{h}
\end{calcu}
\]
The proof of $\psi'\circ \psi\sim \text{id}_{\sum_{x:A}P(x)}<:$ is similar. \\ [0.1cm]
We can use $\Sigma${\sc eq1}, $\Sigma${\sc eq2} and transitivity of equivalence to derive the following inference rule which we will be using later:
\begin{equation}\label{EqProd}
	\inferrule*[right=eq$_{\times}$]
	{f:A\simeq B \\ g: C\simeq D}
	{f\times g:A\times C \simeq B\times D} 
\end{equation} 

\textit{Proof of $+${\sc Eq}1.}Suppose that $\boldsymbol{f}:A\simeq B$. Let  
$\psi:A+C\rightarrow B+C$
be defined by $\psi:\equiv \boldsymbol{\kappa}(\text{inl}\circ f, 
\text{inr}\circ \text{id}_C)$, and let $\psi':B+C\rightarrow A+C$
be defined by $\psi':\equiv \boldsymbol{\kappa}(\text{inl}\circ f', 
\text{inr}\circ\text{id}_C)$.  Let us see that  $\psi'$ is a quasi-inverse of 
$\psi$. Observe that, by definition of $\Psi$ and $\Psi'$, we have 
\begin{multicols}{2}
	$
	\begin{array}{rl}
	\phantom{\equiv}&\psi(\psi'(\text{inl}(x)))\\
	\equiv & \psi(\boldsymbol{\kappa}(\text{inl}\circ f', \text{inr}\circ\text{id}_C)(\text{inl}(x)))\\
	\equiv & \psi(\text{inl}(f'(x)))\\
	\equiv & \boldsymbol{\kappa}(\text{inl}\circ f, \text{inr}\circ \text{id}_C)(\text{inl}(f'(x))\\ 
	\equiv& \text{inl}(f(f'(x))),\quad \text{and}
	\end{array}
	$
	\columnbreak
	\begin{equation}\label{phis}
	\begin{array}{rl}
	\phantom{\equiv}&\psi(\psi'(\text{inr}(y))) \\
	\equiv& \psi(\boldsymbol{\kappa}(\text{inl}\circ f', \text{inr}\circ\text{id}_C)(\text{inr}(y)))\\
	\equiv&\psi(\text{inr}(y))\\
	\equiv& \boldsymbol{\kappa}(\text{inl}\circ f, \text{inr}\circ \text{id}_C)(\text{inr}(y))\\
	\equiv&\text{inr}(y).
	\end{array}
	\end{equation}
\end{multicols}
Then we have
\[
\begin{calcu}
\expro{\psi\circ \psi'\sim \text{id}}
\\
\explo{\equiv}{Definition of $\sim $}
\\
\expro{\prod\limits_{p:B+C}\psi(\psi'(p))=p}
\\
\explo{\leftarrow}{\!:\,$\boldsymbol{\kappa}$}
\\
\expro{\prod\limits_{x:B}(\psi(\psi'(\text{inl}(x)))=\text{inl}(x))\times 
\prod\limits_{y:C}\psi(\psi'(\text{inr}(y)))=\text{inr}(y)}
\\
\explo{\equiv}{Definition of $\psi$ and $\psi'$ (\ref{phis})}
\\
\expro{\prod\limits_{x:B}(\text{inl}(f(f'(x)))=\text{inl}(x))\times 
\prod\limits_{y:C}\text{inr}(y)=\text{inr}(y)}
\\
\explo{\leftarrow}{\!:\,$k$\,\,;\, $k(u,v):\equiv (\lambda x.\text{ap}_{\text{inl}}(u(x)),\lambda x.\text{ap}_{\text{inr}}(v(x))$}
\\
\expro{\prod\limits_{x:B}(f(f'(x))=x)\times \prod\limits_{y:C}y=y}
\\
\explo{\stackrel{\mbox{\tiny $\wedge$}}{\mbox{\tiny :}}}{$h:\equiv (\alpha,\text{refl})$}
\\
\expro{h}
\end{calcu}
\]
We can prove $h'\circ h\sim \text{id}_{\prod_{x:A} P(x)}<:$ similarly.\\ [0.1cm]
\textit{Proof of $+${\sc Eq}2.}
\[
\begin{calcu}
\expro{C+A}
\\
\explo{\simeq}{Commutativity of + (\ref{comm+})}
\\
\expro{A+C}
\\
\explo{\simeq}{$+${\sc Eq}1}
\\
\expro{B+C}
\\
\explo{\simeq}{Commutativity of + (\ref{comm+})}
\\
\expro{C+B}
\end{calcu}
\]
\section{Induction operators as equivalences}
In order to be able to restate HoTT giving equality and equivalence a 
preeminent role, it is convenient (and possible) to show that the inductive 
operators for the equality type, the $\Sigma$-type and the coproduct are 
actually, equivalences. We now proceed to show that this is actually so.

\subsection{Identity type induction operator}\label{IdType}
We prove that for all $P: A\rightarrow {\cal U}$, $\boldsymbol{\iota}$ 
is 
an equivalence, and then,
\[
\prod_{x,y:A}(\prod_{p:x=y}P(x,y,p)) \simeq \prod_{x:A}P(x,x,\text{refl}_x) <:
\]
This equivalence is precisely $\Pi$-[\textbf{Equality}] rule in section \ref{Eind}.

Recall that 
\[
\boldsymbol{\iota}:(\prod_{x:A}P(x,x,\text{refl}_x))\rightarrow
\prod_{x,y:A}\prod_{p:x=y}P(x,y,p).
\]
Now, let us define 
\[
k:\prod_{x,y:A}(\prod_{p:x=y}P(x,y,p))\rightarrow 
\prod_{x:A}P(x,x,\text{refl}_x)
\]
by
\[
k(v)(x):\equiv v(x,x,\text{refl}_x).
\]
Let us prove that $k\circ \boldsymbol{\iota}\sim \text{id}$ and that 
$\boldsymbol{\iota}\circ k\sim \text{id}$.
First, observe that for all $u\!:\!\prod_{x:A}P(x,x,\text{refl}_x)$, by definition of $k$ and $\boldsymbol{\iota}$,
\begin{equation}\label{label1.8.1} 
k(\boldsymbol{\iota}(u))(x)\;\equiv\;\boldsymbol{\iota}(u)(x,x.\text{refl}_x)\;\equiv\; u(x),
\end{equation}
and for all $v\!:\!\prod_{x,y:A}\prod_{p:x=y}P(x,y,p)$,
\begin{equation}\label{label2.8.1}
\boldsymbol{\iota}(k(v))(x,x.\text{refl}_x)\;\equiv\; k(v)(x)\;\equiv\; v(x,x,\text{refl}_x).
\end{equation}
Then, in one hand, because of (\ref{label1.8.1}), we have that $k\circ \boldsymbol{\iota}\sim \text{id}$. On 
the 
other, for each $v:\prod\limits_{x,y:A}\prod\limits_{p:x=y}P(x,y,p)$, let us 
show that $\boldsymbol{\iota}(k(v))=v<:$

\[
\begin{calcu}
\expro{	\boldsymbol{\iota}(k(v))=v}
\\
\explo{\simeq}{Function extensionality (\ref{FuncExt})}
\\
\expro{\prod\limits_{x,y:A}\prod\limits_{p:x=y}\boldsymbol{\iota}(k(v))(x,y,
	p)=v(x,y,p)}
\end{calcu}
\]
\[
\begin{calcu}
\explo{\leftarrow}{\!:\,\, $\boldsymbol{\iota}$}
\\
\expro{\prod\limits_{x:A}\boldsymbol{\iota}(k(v))(x,x,\text{refl}_x)=v(x,x,
	\text{refl}_x)}
\\
\explo{\equiv}{See computation (\ref{label2.8.1}) above}
\\
\expro{\prod\limits_{x:A}v(x,x,\text{refl}_x)=v(x,x,\text{refl}_x)}
\\
\explo{\stackrel{\mbox{\tiny $\wedge$}}{\mbox{\tiny :}}}{$u(x)=\text{refl}_{v(x,x,\text{refl}_x)}$}
\\
\expro{u}
\end{calcu}
\]
Therefore, the equivalence is proved.
\subsection{Identity type based-path induction operator}\label{BasePath}

Let us suppose that $a\!:\!A$ and that $D:\prod_{x:A}\prod_{p:a=x}{\cal U}$. 
Based path induction states the existence of a function $\boldsymbol{\iota}'$ 
presented by the following consequence link
\[
\begin{calcu}
\expro{\prod\limits_{x:A}\prod\limits_{p:a=x}D(x,p)}
\\
\explo{\leftarrow}{\!:\,\,$\boldsymbol{\iota}'_D$;\,\, 
	$\boldsymbol{\iota}'_D(z)(a,\text{refl}_a):\equiv z$}
\\
\expro{D(a,\text{refl}_a)}
\end{calcu}
\]
We have also that $\boldsymbol{\iota}'_P$, the based path induction operator, 
is 
an equivalence, and then
\[\prod_{x:A}(\prod_{p:a=x}P(x,p)) \simeq P(a,\text{refl}_a)<:
\]
This equivalence corresponds to $\Pi$-[\textbf{One-Point}] rule in section \ref{Eind}.\\[0.1cm]
Let us prove that the functions 
\[
\begin{calcu}
\expro{\prod_{x:A}\prod_{p:a=x}P(x,p)}
\\
\explo{\leftarrow}{\!:\,$\boldsymbol{\iota }'\,\,; \,\,\boldsymbol{\iota 
	}'(u)(a,\text{refl}_a):\equiv u$}
\\
\expro{P(a,\text{refl}_a)}
\end{calcu}
\]
and
\[
\begin{calcu}
\expro{P(a,\text{refl}_a)}
\\
\explo{\leftarrow}{\!:\, $k\,\,; \,\,k(v):\equiv v(a,\text{refl}_a)$}
\\
\expro{\prod_{x:A}\prod_{p:a=y}P(x,p)}
\end{calcu}
\]
are quasi-inverses. In fact,
\[k(\boldsymbol{\iota}'(u))\equiv 
\boldsymbol{\iota}'(u)(a,\text{refl}_a)\equiv u,
\]
which shows that 
$k\circ\boldsymbol{\iota}'\sim \text{id}$, and 
\begin{equation}\label{CalBasPath}
\boldsymbol{\iota}'(k(v))(a,\text{refl}_a)\equiv k(v)(x)\equiv 
v(a,\text{refl}_a).
\end{equation}
And so, to prove $\boldsymbol{\iota}'\circ k\sim \text{id}$, it is enough to 
perform the following calculation for all $v:\prod_{x:A}\prod_{p:a=x}P(x,p)$,

\[
\begin{calcu}
\expro{\boldsymbol{\iota}'(k(v))=v}
\\
\explo{\simeq}{Function extensionality (\ref{FuncExt})}
\\
\expro{\prod\limits_{x:A}\prod\limits_{p:a=x}\boldsymbol{\iota}'(k(v))(x,
	p)=v(x,p)}
\end{calcu}
\]
\[
\begin{calcu}
\explo{\leftarrow}{\!: $\boldsymbol{\iota}'$}
\\
\expro{\boldsymbol{\iota}'(k(v))(a,\text{refl}_a)=v(a,\text{refl}_a)}
\\
\explo{\equiv}{See (\ref{CalBasPath}), above}
\\
\expro{v(a,\text{refl}_a)=v(a,\text{refl}_a)}
\\
\explo{\stackrel{\mbox{\tiny $\wedge$}}{\mbox{\tiny :}}}{Definition of $\text{refl}$}
\\
\expro{\text{refl}_{v(a,\text{refl}_a)}}

\end{calcu}
\]
Therefore,
$\prod_{x:A}\prod_{p:a=y}P(x,p)\simeq P(a,\text{refl}_a)<:$

\subsection{$\sum$-type induction operator}

Now, we prove that, for all $P: A\rightarrow {\cal U}$, $\boldsymbol{\sigma}$, the $\sum$-type induction operator, is an equivalence. 
And so,
\begin{equation}\label{SConRule} 
(\prod_{x:A}\prod_{y:B(x)}P((x,y))) \simeq \prod_{g:\sum_{x:A}B(x)}P(g)<: 
\end{equation}
For the case of $P$ being a non-dependent type, the intuitionistic logical 
theorem corresponding to this equivalence is 
\[
\cuan{\forall}{x\!:\!T}{B}{P} \equiv \cuant{\exists}{x\!:\!T}{B} \Rightarrow P  
\]
where $x$ does not occur free in $P$.\\ [0.1cm] 
This motivate us to call the equivalence (\ref{SConRule})  {\it $\Sigma$-consequent rule}.\\[1mm]
Recall that
\[
\boldsymbol{\sigma}:(\prod_{x:A}\prod_{y:B(x)}P((x,y))) 
\rightarrow 
\prod_{g:\sum_{x:A}B(x)}P(g)
\]
and $\boldsymbol{\sigma}(u)((x,y)):\equiv u(x)(y)$. Let 
\[
\Phi:(\prod\limits_{g:\sum_{x:A}B(x)}P(g)) \rightarrow 
\prod\limits_{x:A}\prod\limits_{y:B(x)}P((x,y))
\]
be defined by $\Phi(v)(x)(y):\equiv v((x,y))$. Composing 
$\boldsymbol{\sigma}$ with $\Phi$ we get
\[
\Phi(\boldsymbol{\sigma}(u))(x)(y)\equiv \boldsymbol{\sigma}(u)((x,y))\equiv u(x)(y).
\]
Then $\Phi\circ \boldsymbol{\sigma}$ is homotopic to the identity function. 
Conversely, let $v$ be an inhabitant of $\prod_{g:\sum_{x:A}B(x)}P(g)$, then
\[
\begin{calcu}
\expro{\boldsymbol{\sigma}(\Phi (v))=v}
\\
\explo{\simeq}{Function extensionality (\ref{FuncExt})}
\\
\expro{\prod\limits_{g:\sum_{x:A}B(x)}\boldsymbol{\sigma}(\Phi(v))(g)=v(g)}
\\
\explo{\leftarrow}{\!:\,$\boldsymbol{\sigma}$}
\\
\expro{\prod\limits_{x:A}\prod\limits_{y:B(x)}\boldsymbol{\sigma}(\Phi(v))(x,y)=v((x,y))}
\\
\explo{\equiv}{$\boldsymbol{\sigma}(\Phi(v))((x,y))\equiv \Phi(v)(x)(y) \equiv v((x,y))$}
\\
\expro{\prod\limits_{x:A}\prod\limits_{y:B(x)}v((x,y))=v((x,y))}
\\
\explo{\stackrel{\mbox{\tiny $\wedge$}}{\mbox{\tiny :}}}{$h:\equiv \lambda x.\lambda y.\text{refl}_v(x,y)$}
\\
\expro{h}
\end{calcu}
\]
So, $\boldsymbol{\sigma}\circ\Phi$ is homotopic to the identity function. 

\subsection{Coproduct induction operator}\label{CopIndOp}

For all $A,B:\mathcal U$ and $P:A+B\rightarrow 
\mathcal U$ we have that 
\[
(\prod\limits_{x:A+B}P(x)) \simeq 
(\prod\limits_{x:A}P(\text{inl}(x)))\times 
\prod\limits_{y:B}P(\text{inr}(x))<:
\]
This equivalence correspond to $\Pi$-[\textbf{Range Split}] rule in 
section \ref{Eind}.\\[0.5mm]

\textit{Proof.} We have the induction operator $\boldsymbol{\kappa}$:
\[
\begin{calcu}
\expro{\prod\limits_{x:A+B}P(x)}
\\
\explo{\leftarrow}{\!:\,$\boldsymbol{\kappa}$\,;\, $\boldsymbol{\kappa}(u,v)(\text{inl}(x)):\equiv u(x)\,;\,\boldsymbol{\kappa}(u,v)(\text{inr}(x)):\equiv v(x)$}
\\
\expro{(\prod\limits_{x:A}P(\text{inl}(x))\times 
\prod\limits_{y:B}P(\text{inr}(x))}
\end{calcu}
\]
and let us define 
\[
\Psi: (\prod\limits_{x:A+B}P(x))\,\to\, 
(\prod\limits_{x:A}P(\text{inl}(x)))\times 
\prod\limits_{y:B}P(\text{inr}(y))
\]
by $\Psi(g):\equiv (g\circ \text{inl},g\circ \text{inr})$. 
Let us see that  $\Psi$ is a quasi-inverse of $\boldsymbol{\kappa}$. We 
show that, the type $\boldsymbol{\kappa}\circ\Psi\sim \text{id}$, which by 
definition is equivalent to 
\[
\prod\limits_{g:\prod\limits_{x:A+B}P(x)}\boldsymbol{\kappa}(\Psi(g)))=g,
\]
is inhabited. Let $g$ be an object of type $\prod_{x:A+B}P(x)$, then:
\[
\begin{calcu}
\expro{\boldsymbol{\kappa}(\Psi(g)))=g}
\\
\explo{\equiv}{Definition of $\Psi$}
\\
\expro{\boldsymbol{\kappa}(g\circ \text{inl}, g\circ \text{inr})=g}
\\
\explo{\simeq}{Function extensionality (\ref{FuncExt})}
\\
\expro{\boldsymbol{\kappa}(g\circ \text{inl}, g\circ \text{inr})\sim g}
\\
\explo{\equiv}{Definition of $\sim$}
\\
\expro{\prod\limits_{z:A+B}\boldsymbol{\kappa}(g\circ \text{inl}, g\circ 
\text{inr})(z
)=g(z)}
\\
\explo{\leftarrow}{\!:\,$\boldsymbol{\kappa}$}
\\
\expro{\phantom{\times}\prod\limits_{x:A}\boldsymbol{\kappa}(g\circ 
\text{inl}, g\circ \text{inr})(\text{inl} (x))=g(\text{inl} (x))\\ 
&\times\prod\limits_{y:B}\boldsymbol{\kappa}(g\circ \text{inl}, g\circ 
\text{inr})(\text{inr}(y))=g(\text{inr}(y))}
\\
\explo{\equiv}{Definition of $\boldsymbol{\kappa}$}
\\
\expro{\prod\limits_{x:A}((g\circ \text{inl})(x)=( g\circ \text{inl})(x)) 
\times \prod\limits_{y:B}(g\circ \text{inr})(y)=( g\circ \text{inr})(y)}\\
\explo{\stackrel{\mbox{\tiny $\wedge$}}{\mbox{\tiny :}}}{$h(g):\equiv (\lambda x.\text{refl}_{g(\text{inl}(x))},\lambda x.\text{refl}_{g(\text{inr}(t))})$}
\\
\expro{h(g)}
\\
\end{calcu}
\]
And now, we show that $\Psi\circ\boldsymbol{\kappa}\sim 
\text{id}<:$\,. In other words, that
\[
\prod\limits_{u:\prod\limits_{x:A}P(\text{inl}(x))\times 
\prod\limits_{y:B}P(\text{inr}(y))}\Psi(\boldsymbol{\kappa}(u))=u\,<:
\]
Let $u$ be an object of type $\prod_{x:A}P(\text{inl}(x))\times 
\prod_{y:B}P(\text{inr}(y))$, $p:u\!=\!(\text{pr}_1(u),\text{pr}_2(u))$ and 
 $Q$ the type family defined by $Q(u):\equiv (\boldsymbol{\kappa}(u)\!\circ\! 
\text{inl},\,\boldsymbol{\kappa}(u) \!\circ\! 
\text{inr})\!=\!u $, and so, by the second Leibniz principle,
\[
\text{tr}^Q(u,(\text{pr}_1(u),\text{pr}_2(u)),p) : Q(u)\simeq Q((\text{pr}_1(u),\text{pr}_2(u)))
\]
Then:
\[
\begin{calcu}
\expro{\Psi(\boldsymbol{\kappa}(u))=u}
\\
\explo{\equiv}{Definition of $\Psi$}
\\
\expro{(\boldsymbol{\kappa}(u)\circ \text{inl},\boldsymbol{\kappa}(u)\circ 
\text{inr})=u}
\\
\explo{\simeq}{\!:\, $\text{tr}^Q(u,(\text{pr}_1(u),\text{pr}_2(u)),p)$}
\\
\expro{(\boldsymbol{\kappa}(\text{pr}_1(u),\text{pr}_2(u))\circ 
\text{inl},\boldsymbol{\kappa}(\text{pr}_1(u),\text{pr}_2(u))\circ \text{inr})=(\text{pr}_1(u),\text{pr}_2(u))}
\\
\explo{\simeq}{$(a,b)=(c,d)\simeq (a=c)\times (b=d)\,<:$}
\\
\expro{(\boldsymbol{\kappa}(\text{pr}_1(u),\text{pr}_2(u))\circ \text{inl}=\text{pr}_1(u)) \times 
(\boldsymbol{\kappa}(\text{pr}_1(u),\text{pr}_2(u))\circ \text{inr}=\text{pr}_2(u))}
\\
\explo{\equiv}{Definition of $\boldsymbol{\kappa}$}
\\
\expro{(\text{pr}_1(u)=\text{pr}_1(u))\,\,\times\,\, (\text{pr}_2(u)=\text{pr}_2(u))}
\\
\explo{\stackrel{\mbox{\tiny $\wedge$}}{\mbox{\tiny :}}}{$h:\equiv \text{refl}_{\text{pr}_1(u)}\,\,;\,\,k:\equiv \text{refl}_{\text{pr}_2(u)}$}
\\
\expro{(h,k)}
\\
\end{calcu}
\] 
As a matter of fact, the induction operators corresponding to $W$ type, $\mathds O$ type and  $\mathds 1$ type could be similarly proved to be equivalences. 
\section{Operational properties of $\Pi$ and $\Sigma$ types} 
Now we come back to the operational rules enumerated in section \ref{Eind} and prove the ones that we have not proved yet.\\ [0.1cm] 
[\textbf{One-Point}] rules. In first order logic, quantifying a property over exactly one element is equivalent to the property applied to just this element. For the case of HoTT, this properties are slightly more general.
\[\prod_{x:A}(\prod_{p:a=x}P(x,p)) \simeq P(a,\text{refl}_a)<:
\]
and
\[
\sum_{x:A}(\sum_{p:x=a}P(x,p))\simeq P(a,\text{refl}_a)<:.
\]
We have proved $\Pi$-[\textbf{One-Point}] rule in subsection \ref{BasePath}. We now prove 
the $\Sigma$-[\textbf{One-Point}] rule.\\ [0.1cm]

Given $A:{\cal U}$, $a:A$ and $P:\prod_{x:A}\prod_{p:x=a}{\cal U}$, let us 
construct
\[
\Phi:\sum_{x:A}(\sum_{p:x=a}P(x,p))\rightarrow P(a,\text{refl}_a).
\]
This can be done by means of the following deductive chain:
\[
\begin{calcu}
\expro{\prod\limits_{g:\sum_{x:A}\sum_{p:x=a}P(x,p)} P(a,\text{refl}_a)}
\\
\explo{\simeq}{\!:\,$\boldsymbol{\sigma}$,\, $\Sigma$-consequent rule}
\\
\expro{\prod\limits_{x:A}\prod\limits_{y:\sum_{p:x=a}P(x,p)}P(a,\text{refl}_a)}
\\
\explo{\simeq}{\!:\,$\Delta\,;\,\varphi_x:\equiv \boldsymbol{\sigma}_{x}$,\, $\Pi${\sc eq1}}
\\
\expro{\prod\limits_{x:A}\prod\limits_{p:x=a}\prod\limits_{z:P(x,p)}P(a,\text{
		refl}_a)}
\\
\explo{\simeq}{\!:\,$\boldsymbol{\iota}'$,\, $\Pi$-one-point rule}
\\
\expro{\prod\limits_{z:P(a,\text{refl}_a)}P(a,\text{refl}_a)}
\\
\explo{\stackrel{\mbox{\tiny $\wedge$}}{\mbox{\tiny :}}}{$u:\equiv \text{id}_{P(a,\text{refl}_a)}$}
\\
\expro{u}
\end{calcu}
\]
In the chain above, $\boldsymbol{\sigma}_{x}$ is the induction operator for $\sum_{p:x=a}P(x,p)$ evaluated at the constant type 
family $C(x,y):\equiv P(a,\text{refl}_a)$. \\ [0.1cm]
Now, let $\Psi:P(a,\text{refl}_a)\rightarrow \sum_{x:A}\sum_{p:x=a}P(x,p) $ be 
defined by  \[\Psi(u):\equiv (a,(\text{refl}_a,u)).\]
Let us verify that $\Phi\circ\Psi\sim \text{id}$ and that $\Psi\circ\Phi\sim 
\text{id}$. First of all observe that,
making the compositions in the above chain, we get  
\[
\Phi :\equiv \boldsymbol{\sigma}(\Delta(\boldsymbol{\iota}'(\text{id}_{P(a,\text{refl}_a)}))).
\]
On one hand we have,
\begin{multicols}{2}
	$\begin{array}{rl}
	\phantom{\equiv}& \Phi(\Psi(t))\\
	\equiv &\boldsymbol{\sigma}(\Delta(\boldsymbol{\iota}'(\text{id}_{P(a,\text{refl}_a)}))(a,\text{refl}_a,t))\\
	\equiv &\Delta(\boldsymbol{\iota}'(\text{id}_{P(a,\text{refl}_a)}))(a)((\text{refl}_a,t))\\
	\end{array}$
	\columnbreak
	$\begin{array}{rll}
	\equiv &\boldsymbol{\sigma}_a(\boldsymbol{\iota}'(\text{id}_{P(a,\text{refl}_a)})(a))((\text{refl}_a,t))\\
	\equiv& \boldsymbol{\iota}'(\text{id}_{P(a,\text{refl}_a)})(a)(\text{refl}_a)(t)\\
	\equiv& \text{id}_{P(a,\text{refl}_a)}(t)
	\equiv t&
	\end{array}$
\end{multicols}
and, on the other hand, 

\[
\begin{calcu}
\expro{\Psi\circ\Phi \sim \text{id}}
\\
\explo{\equiv}{Definition of $\sim$}
\\
\expro{\prod\limits_{g:\sum_{x:A}\sum_{p:x=a}P(x,p)} \Psi(\Phi(g))=g}
\\
\explo{\equiv}{Definition of $\Psi$}
\\
\expro{\prod\limits_{g:\sum_{x:A}\sum_{p:x=a}P(x,p)} 
	(a,(\text{refl}_a,\Phi(g))=g}
\\
\explo{\simeq}{\!:\,$\boldsymbol{\sigma}$,\, $\Sigma$-consequent rule}
\\
\expro{\prod\limits_{x:A}\prod\limits_{y:\sum_{p:x=a}P(x,p)}(a,(\text{refl}_a,
	\Phi((x,y)))=(x,y)}
\end{calcu}
\]
\[
\begin{calcu}
\explo{\simeq}{\!:\,$\Delta$\,\,;\,\,$\varphi_x:\equiv \boldsymbol{\sigma}_{x}$\,;\, $\Pi${\sc eq1}}
\\
\expro{\prod\limits_{x:A}\prod\limits_{p:x=a}\prod\limits_{z:P(x,p)}(a,(\text
{refl}_a,\Phi((x,(p,z))))=(x,(p,z))}
\\
\explo{\simeq}{\!:\,$\boldsymbol{\iota}'$ ($\Pi$-one-point rule)}
\\
\expro{\prod\limits_{z:P(a,\text{refl}_a)}(a,(\text{refl}_a,\Phi((a,(\text{refl}
	_a,z))))=(a,(\text{refl}_a,z))}
\\
\explo{\equiv}{Property of $\Phi$}
\\
\expro{ \prod\limits_{z: 
		P(a,\text{refl}_a)}(a,(\text{refl}_a,z))=(a,(\text{refl}_a,z))}
\\
\explo{\stackrel{\mbox{\tiny $\wedge$}}{\mbox{\tiny :}}}{$h(z):\equiv \text{refl}_{(a,(\text{refl}_a,z))}$}
\\
\expro{h}
\end{calcu}
\]

[\textbf{Equality}] rules. These equivalences correspond, in first 
order logic, to the case when we are quantifying over two variables that happen 
to be equal, then one of those quantified variables may be made equal to the 
other, and be, in this way, eliminated.  
\[
\prod_{x,y:A}(\prod_{p:x=y}P(x,y,p)) \simeq \prod_{x:A}P(x,x,\text{refl}_x)<:
\]
and
\[
\sum_{x,y:A}(\sum_{p:x=y}P(x,y,p)) \simeq \sum_{x:A}P(x,x,\text{refl}_x)<:
\]
$\Pi$-[\textbf{Equality}] rule was proved in subsection \ref{IdType}. The proof of 
$\Sigma$-[\textbf{Equality}] rule follows analogous steps to those of the 
$\Sigma$-[\textbf{One-Point}] rule. We omit it.\\ [0.1cm]
[\textbf{Range Split}] rules. The range split rule is a property of 
operationals in general. In the case of logical quantifications, it allows separating them into two quantifiers of the same kind of the original one: universal or existential. These operational parts are  joined by conjunctions for the first kind, and by disjunctions for the second. Their ranges correspond to disjoint components of the range of the original quantification. In the case of HoTT, this splitting is possible when the range of a $\Pi$-type or a $\Sigma$-type corresponds to a coproduct type. For the case 
of a $\Pi$-type, $\Pi$-[\textbf{Range Split}], its parts are joined by a Cartesian product and in the case of a 
$\Sigma$-type, $\Sigma$-[\textbf{Range Split}], they are joined by a coproduct operator, namely,
	\[
	\prod\limits_{x:P+Q}R(x)\simeq (\prod\limits_{x:P}R(\text{inl}(x)) 
) 
\times (\prod\limits_{x:Q}R(\text{inr}(x)) )
	\]
	and
	\[\sum\limits_{x:P+Q}R(x)\simeq (\sum\limits_{x:P}R(\text{inl}(x)) ) 
+(\sum\limits_{x:Q}R(\text{inr}(x)) )
	\]
The $\Pi$-[\textbf{Range Split}] rule is related to the coproduct induction operator and  was proved in subsection \ref{CopIndOp}. We now prove $\Sigma$-[\textbf{Range Split}] rule.\\ [0.1cm]
In order to get a function
$$\Phi:(\sum_{x:P+Q}R(x))\,\rightarrow\, 
(\sum_{y:P}R(\text{inl}(y)))+\sum_{z:Q}R(\text{inr}(z))
$$
let us consider the following deductive chain:
\[
\begin{calcu}
\expro{(\sum\limits_{x:P+Q}R(x))\rightarrow 
	(\sum\limits_{y:P}R(\text{inl}(y)))+\sum\limits_{z:Q}R(\text{inr}(z))}
\\
\explo{\simeq}{\!:\,$\boldsymbol{\sigma}$,\, $\Sigma$-consequent rule}
\\
\expro{\prod\limits_{x:P+Q}(R(x) \,\rightarrow\, 
	(\sum\limits_{y:P}R(\text{inl}(y))+\sum\limits_{z:Q}R(\text{inr}(z))))}
\end{calcu}
\]
\[
\begin{calcu}
\explo{\simeq}{\!:\,$\boldsymbol{\kappa}$, ($\Pi$-range split rule)}
\\
\expro{\phantom{\times} (\prod\limits_{u:P}(R(\text{inl}(u))\rightarrow 
	(\sum\limits_{y:P}R(\text{inl}(y)))+\sum\limits_{z:Q}R(\text{inr}(z)))) }
& \times (\prod\limits_{v:Q}(R(\text{inr}(v))\rightarrow 
(\sum\limits_{y:P}R(\text{inl}(y)))+\sum\limits_{z:Q}R(\text{inr}(z)))) 
\\
\\
\explo{\stackrel{\mbox{\tiny $\wedge$}}{\mbox{\tiny :}}}{$\phi_0(u)(a):\equiv \text{inl}((u,a));\quad 
	\phi_1(v)(b):\equiv \text{inr}((v,b))$}
\\
\expro{(\phi_0,\phi_1)}
\end{calcu}
\]
Then we can put $\Phi:\equiv 
\boldsymbol{\sigma}(\boldsymbol{\kappa}(\phi_0,\phi_1))$\\ [0.1cm]
Now, in order to get a function  
\[\Psi:\sum_{y:P}R(\text{inl}(y))+\sum_{z:Q}R(\text{inr}(z))\rightarrow 
\sum_{x:P+Q}R(x)
\]
let us consider the following deductive chain:
\[
\begin{calcu}
\expro{(\sum\limits_{y:P}R(\text{inl}(y))+(\sum\limits_{z:Q}R(\text{inr}
	(z))\rightarrow \sum\limits_{x:P+Q}R(x)}
\\
\explo{\simeq}{\!:\,$\boldsymbol{\kappa}$,\, ($\Pi$-range split rule)}
\\
\expro{((\sum\limits_{y:P}R(\text{inl}(y))\rightarrow 
	\sum\limits_{x:P+Q}R(x)) 
	\,\times\, 
	((\sum\limits_{z:Q}R(\text{inr}(z))\rightarrow 
	\sum\limits_{x:P+Q}R(x)) }
\\
\explo{\simeq}{\!:\,$\boldsymbol{\sigma}_1\!\times\!\boldsymbol{\sigma}_2$,\, {\sc eq}$_{\times}$ (\ref{EqProd})}
\\
\expro{((\prod\limits_{y:P}R(\text{inl}(y))\rightarrow 
	\sum\limits_{x:P+Q}R(x)) \times ((\prod\limits_{z:Q}
	R(\text{inr}(z))\rightarrow \sum\limits_{x:P+Q}R(x)) }
\\
\explo{\stackrel{\mbox{\tiny $\wedge$}}{\mbox{\tiny :}}}{$\psi_0(y)(a):\equiv (\text{inl}(y),a);\quad 
	\psi_1(z)(b):\equiv (\text{inr}(z),b)$}
\\
\expro{(\psi_0,\psi_1)}
\end{calcu}
\]
Then we may define $\Psi:\equiv \boldsymbol{\kappa}(\boldsymbol{\sigma}_1\!\times\!\boldsymbol{\sigma}_2(\psi_0,\psi_1)):\equiv \boldsymbol{\kappa}(\boldsymbol{\sigma}_1(\psi_0),\boldsymbol{\sigma}_2(\psi_1))$\\ [0.1cm]
Observe that 
\begin{multicols}{2}
	$\begin{array}{rl}	
	\phantom{\equiv} &\Phi(\Psi(\text{inl}(f_1,f_2)))\\
	\equiv&\Phi(\boldsymbol{\kappa}(\boldsymbol{\sigma}_1(\psi_0),\boldsymbol{\sigma}_2(\psi_1))(\text{inl}(f_1,f_2)))\\
	\equiv&\Phi(\boldsymbol{\sigma}_1(\psi_0)(f_1,f_2))\\
	\equiv&\Phi(\psi_0(f_1)(f_2))\\
	\end{array}$
	
	\columnbreak
	
	$\begin{array}{rl}	
	\equiv&\Phi(\text{inl}(f_1),f_2)\\
	\equiv&\boldsymbol{\kappa}(\phi_0,\phi_1))(\text{inl}(f_1))(f_2)\\
	\equiv&\phi_0(f_1)(f_2)\\
	\equiv&\text{inl}(f_1,f_2).
	\end{array}$
\end{multicols}
%
%
%
%
In the same way we can prove that $\Phi(\Psi(\text{inr}(g_1,g_2)))\equiv \text{inr}(g_1,g_2))$
Then

\[
\begin{calcu}
\expro{\prod\limits_{p:\sum_{y:P}R(\text{inl}(y))+\sum_{z:Q}R(\text{inr}(z))}
	\Phi(\Psi(p))=p}
\\
\explo{\simeq}{\!:\,$\boldsymbol{\kappa}$,\, ($\Pi$-range split rule)}
\\
\expro{\phantom{\times}\prod\limits_{f:\sum_{y:P}R(\text{inl}(y))}
	\Phi(\Psi(\text{inl}(f)))=\text{inl}(f)\\
	&\times 
	\prod\limits_{g:\sum_{x:Q}R(\text{inr}(x))}\Phi(\Psi(\text{inr}(g)))=\text{inr}
	(g)}
\\
\explo{\simeq}{\!:\, $\boldsymbol{\sigma}_1\times \boldsymbol{\sigma}_2$, {\sc eq}$_{\times}$ (\ref{EqProd})}
\\
\expro{\phantom{\times}\prod\limits_{f_1:P}\prod\limits_{f_2:R(\text{inl}(f_1))}
	\Phi(\Psi(\text{inl}(f_1,f_2)))=\text{inl}(f_1,f_2)\\
	&\times
		\prod\limits_{g_1:P}\prod\limits_{g_2:R(\text{inr}(g_1))}\Phi(\Psi(\text{inr}
	(g_1,g_2)))=\text{inr}(g_1,g_2)}
\end{calcu}
\]
\[
\begin{calcu}
\explo{\equiv}{Above computations}
\\
\expro{\phantom{\times}\prod\limits_{f_1:P}\prod\limits_{f_2:R(\text{inl}(f_1))}
	\text{inl}(f_1,f_2)=\text{inl}(f_1,f_2)\\
	&\times 
	\prod\limits_{g_1:P}\prod\limits_{g_2:R(\text{inr}(g_1))}\text{inr}(g_1,
	g_2)=\text{inr}(g_1,g_2)}
\\
\explo{\stackrel{\mbox{\tiny $\wedge$}}{\mbox{\tiny :}}}{$u(f_1,f_2):\equiv 
	\text{refl}_{\textbf{inl}(f_1,f_2)}\,\,;\,\,u(g_1,g_2):\equiv 
	\text{refl}_{\textbf{inr}(g_1,g_2)}$}
\\
\expro{(u,v)}
\end{calcu}
\]
In the other direction, observe that
\begin{multicols}{2}
	$\begin{array}{rl}
	\phantom{\equiv}&\Psi(\Phi(\text{inl}(w), u_2))\\
	\equiv&\Psi(\boldsymbol{\sigma}(\boldsymbol{\kappa}(\phi_0,\phi_1))(\text{inl}(w),u_2))\\
	\equiv&\Psi(\boldsymbol{\kappa}(\phi_0,\phi_1)(\text{inl}(w))(u_2))\\
	\equiv&\Psi(\phi_0(w)(u_2))\\
	\end{array}$
	
	\columnbreak
	
	$\begin{array}{rll}
	\equiv&\Psi(\text{inl}(w,u_2))\\
	\equiv&\boldsymbol{\kappa}(\boldsymbol{\sigma}_1(\psi_0),\boldsymbol{\sigma}_2(\psi_1))(\text{inl}(w,u_2))\\
	\equiv&\boldsymbol{\sigma}_1(\psi_0)(w,u_2)\\
	\equiv&\psi_0(w)(u_2)\equiv (\text{inl}(w),u_2).
	\end{array}$
\end{multicols}
In the same way we can prove that $\Psi(\Phi(\text{inr}(z), u_2)):\equiv 
(\text{inr}(z), u_2)$.
Then
\[
\begin{calcu}
\expro{\prod\limits_{u:\sum_{x:P+Q}R(x)}\Psi(\Phi(u))=u}
\\
\explo{\simeq}{\!:\,$\boldsymbol{\sigma}$,\, $\Sigma$-consequent rule}
\\
\expro{\prod\limits_{u_1:P+Q}\prod\limits_{u_2:R(u_1)}\Psi(\Phi(u_1,u_2))=(u_1,
	u_2)}
\\
\explo{\simeq}{\!:\,$\boldsymbol{\kappa}$,\, ($\Pi$-range split rule)}
\\
\expro{\phantom{\times}\prod\limits_{w:P}\prod\limits_{u_2:R(\text{inl}(w))}
	\Psi(\Phi(\text{inl}(w),u_2))=
	(\text{inl}(w),u_2)\\
	&\times 
	\prod\limits_{z:Q}\prod\limits_{u_2:R(\text{inr}(z))}\Psi(\Phi(\text{inr}(z),
	u_2))=(\text{inr}(z),u_2)}
\\
\explo{\equiv}{Above computations}
\\
\expro{\phantom{\times}\prod\limits_{w:P}\prod\limits_{u_2:R(\text{inl}(w))}
	(\text{inl}(w),u_2)=
	(\text{inl}(w),u_2)\\
	&\times 
	\prod\limits_{z:Q}\prod\limits_{u_2:R(\text{inr}(z))}(\text{inr}(z),u_2)=(\text{
		inr}(z),u_2)}
\\
\explo{\stackrel{\mbox{\tiny $\wedge$}}{\mbox{\tiny :}}}{$h:\equiv (\lambda w.\lambda u_2.\text{refl}_{(\text{inl}(w),u_2)}, \lambda z.\lambda u_2.\text{refl}_{(\text{inr}(z),u_2)})$}
\\
\expro{h}
\end{calcu}
\]
[\textbf{Term Split}] rules. In logic, universal quantifications of conjuntions 
split (through an equivalence) into universal quantifications of each conjunct  
joined by conjunctios too. Dually, existential quantifications split into existential quantifications of each disjunct joined by disjunctions. In the case of HoTT, $\Pi$-types mapping into Cartesian products split into $\Pi$-types for each factor joined by Cartesian products, $\Pi$-[\textbf{Term Split}] rule. Dually, for 
$\Sigma$-types, we have an analogous situation replacing cross products by coproducts, $\Sigma$-[\textbf{Term Split}] rule. Namely,
	\[
	\prod\limits_{x:A}(P(x)\times Q(x)) \simeq (\prod\limits_{x:A}P(x) 
) \times (\prod\limits_{x:A}Q(x) )
	\]
	and
	\[
	\sum\limits_{x:A}(P(x)+ Q(x)) \simeq (\sum\limits_{x:A}P(x) 
) 
+ (\sum\limits_{x:A}Q(x) )
	\]	
\\[0.1cm]
To prove $\Pi$-[\textbf{Term Split}] rule, let $\Phi:(\prod_{x:A}P(x)) \times 
(\prod_{y:A}Q(y)) \rightarrow \prod_{x:A}P(x) \times Q(x)$ be 
defined by $\Phi(u)(x)\!:\equiv\! ((\text{pr}_1u)(x),(\text{pr}_2u)(x))$, and also, let
$\Psi\!:\!\prod_{x:A}(P(x)\!\times\! Q(x)) \rightarrow 
(\prod\limits_{x:A}P(x)) \times\prod_{y:A}Q(y)$ be defined by 
$\Psi(g):\equiv (\text{pr}_1\circ g, \text{pr}_2\circ g)$.  Let us see that 
$\Psi$ is a quasi-inverse of $\Phi$:
\[
\begin{calcu}
\expro{\Psi\circ \Phi \sim 
	\text{id}_{\prod\limits_{x:A}P(x)\times\prod\limits_{y:A}Q(y)}}
\\
\explo{\equiv}{Definition of $\sim$}
\\
\expro{\prod\limits_{u:\prod\limits_{x:A}P(x)\times\prod\limits_{y:A}Q(y)}
	\Psi(\Phi(u))=u}
\\
\explo{\equiv}{Definition of $\Psi$}
\\
\expro{\prod\limits_{u:\prod\limits_{x:A}P(x)\times\prod\limits_{y:A}Q(y)}(\text
	{pr}_1 \circ \Phi(u), \text{pr}_2 \circ \Phi(u))=u}
\\
\explo{\equiv}{Definition of $\Phi$}
\\
\expro{\prod\limits_{u:\prod\limits_{x:A}P(x)\times\prod\limits_{y:A}Q(y)}(\text
	{pr}_1u, \text{pr}_2u)=u}
\\
\explo{\stackrel{\mbox{\tiny $\wedge$}}{\mbox{\tiny :}}}{Uniqueness principle of pairs (\ref{UniqPairs})}
\\
\expro{h}
\\
\end{calcu}
\]
Now let us show that $\Phi\circ \Psi \sim \text{id}<:$
\[
\begin{calcu}
\expro{\Phi\circ \Psi \sim \text{id}_{\prod\limits_{x:A}P(x)\times Q(x)}}
\\
\explo{\equiv}{Definition of $\sim$}
\\
\expro{\prod\limits_{g:\prod\limits_{x:A}P(x)\times Q(x)}\Phi(\Psi(g)))=g}
\\
\explo{\equiv}{Definition of $\Psi$}
\\
\expro{\prod\limits_{g:\prod\limits_{x:A}P(x)\times Q(x)}\Phi((\text{pr}_1 
	\circ 
	g, \text{pr}_2 \circ g))=g}
\\
\explo{\simeq}{Function extensionality (\ref{FuncExt})}
\\
\expro{\prod\limits_{g:\prod\limits_{x:A}P(x)\times Q(x)}\Phi((\text{pr}_1 
	\circ 
	g, \text{pr}_2 \circ g))\sim g}
\\
\explo{\equiv}{Definition of $\sim$}
\\
\expro{\prod\limits_{g:\prod\limits_{x:A}P(x)\times 
		Q(x)}\prod\limits_{x:A}\Phi((\text{pr}_1 \circ g, \text{pr}_2 \circ g))(x)= 
	g(x)}
\\
\explo{\equiv}{Definition of $\Phi$}
\\
\expro{\prod\limits_{g:\prod\limits_{x:A}P(x)\times 
		Q(x)}\prod\limits_{x:A}(\text{pr}_1 (g(x)), \text{pr}_2 (g(x))=g(x)}
\\
\explo{\stackrel{\mbox{\tiny $\wedge$}}{\mbox{\tiny :}}}{Uniqueness principle of pairs (\ref{UniqPairs})}
\\
\expro{k.}
\\
\end{calcu}
\]
And now, we prove the $\Sigma$-[\textbf{Term Split}] rule:\\ [0.1cm]
In order to get a function 
\[
\Phi:\sum_{x:A} (P(x)+Q(x)) \rightarrow (\sum_{x:A}P(x)) + (\sum_{x:A}Q(x))
\]
let us consider the folowing deductive chain:
\[
\begin{calcu}
\expro{\sum\limits_{x:A}(P(x)+Q(x)) \rightarrow (\sum\limits_{x:A}P(x)) 
	+(\sum\limits_{x:A}Q(x))}
\\
\explo{\simeq}{\!:\,$\boldsymbol{\sigma}$,\, $\Sigma$-consequent rule}
\\
\expro{\prod\limits_{x:A}((P(x)+Q(x)) \rightarrow (\sum\limits_{x:A}P(x)) 
	+(\sum\limits_{x:A}Q(x)))}
\\
\explo{\simeq}{\!:\,$\Delta\,\,;\,\,\varphi_x:\equiv \boldsymbol{\kappa}_x$, $\Pi$\sc{eq1}}
\\
\expro{\prod\limits_{x:A}((P(x) \rightarrow 
	(\sum\limits_{x:A}P(x)) +(\sum\limits_{x:A}Q(x))) \times (Q(x) \rightarrow 
	(\sum\limits_{x:A}P(x)) +(\sum\limits_{x:A}Q(x))))}
\\
\explo{\simeq}{\!:\,$\eta\,\,;\,\,\eta(u,v):\equiv\lambda x.(u(x),v(x))$, ($\Pi$-term split rule) }
\\
\expro{\prod\limits_{x:A}(P(x) \rightarrow 
	\sum\limits_{x:A}P(x)+\sum\limits_{x:A}Q(x)\times \prod\limits_{x:A}(Q(x) \rightarrow 
	\sum\limits_{x:A}P(x)+\sum\limits_{x:A}Q(x))}
\\
\explo{\stackrel{\mbox{\tiny $\wedge$}}{\mbox{\tiny :}}}{$\phi_1:\equiv\lambda x.\lambda y.\text{inl}(x,y)$;  
	$\phi_2:\equiv\lambda x.\lambda z.\text{inr}(x,z)$}
\\
\expro{(\phi_1,\phi_2)}
\end{calcu}
\]
In the chain above, $\boldsymbol{\kappa}_{x}$ is the induction operator for $P(x)+Q(x)$ evaluated at the constant type 
family $D:(\sum_{x:A}P(x))+(\sum_{x:A}Q(x))$. Then, we may define $\Phi:\equiv 
\boldsymbol{\sigma}(\Delta(\eta(\phi_1,\phi_2)))$.\\ [0.1cm]
In order to get a function 
\[
\Psi:\sum_{x:A}P(x)+\sum_{x:A}Q(x)\rightarrow \sum_{x:A}P(x)+Q(x)
\]
let us consider the following deductive chain: 
\[
\begin{calcu}
\expro{(\sum\limits_{x:A}P(x)) +(\sum\limits_{x:A}Q(x)) \rightarrow 
	\sum\limits_{x:A}(P(x)+Q(x))}
\\
\explo{\simeq}{\!:\,$\boldsymbol{\kappa}$, ($\Pi$-range split rule)}
\\
\expro{((\sum\limits_{x:A}P(x)) \rightarrow \sum\limits_{x:A}P(x)+Q(x)\times 
	((\sum\limits_{x:A}Q(x)) \rightarrow \sum\limits_{x:A}P(x)+Q(x))}
\\
\explo{\simeq}{\!:\, $\boldsymbol{\sigma}_1\!\times\!\boldsymbol{\sigma}_2$, {\sc eq}$_{\times}$ (\ref{EqProd})}
\\
\expro{\prod\limits_{x:A}(P(x) \rightarrow \sum\limits_{x:A}P(x)+Q(x))\times 
	\prod\limits_{x:A}(Q(x) \rightarrow \sum\limits_{x:A}P(x)+Q(x))}
\\
\explo{\stackrel{\mbox{\tiny $\wedge$}}{\mbox{\tiny :}}}{$\psi_1:\equiv\lambda x.\lambda y.(x,\text{inl}(y))$\,;\,  
	$\psi_2:\equiv \lambda x.\lambda z.(x,\text{inr}(z))$}
\\
\expro{(\psi_1,\psi_2)}
\end{calcu}
\]
Then we may define $\Psi:\equiv \boldsymbol{\kappa}(\boldsymbol{\sigma}_1\!\times\!\boldsymbol{\sigma}_2(\psi_1,\psi_2)))$\\ 
Observe that
\begin{multicols}{2}
	$\begin{array}{rl}
	\phantom{\equiv}&\Phi(\Psi(\text{inl}(a_1,a_2)))\\
	\equiv&\Phi(\boldsymbol{\kappa}(\boldsymbol{\sigma}_1(\psi_1),\boldsymbol{\sigma}_2(\psi_2))\text{inl}(a_1,a_2))\\
	\equiv & \Phi(\boldsymbol{\sigma}_1(\psi_1) (a_1,a_2)))\\
	\equiv&\Phi(\psi_1(a_1)(a_2))\\
	\equiv& \Phi(a_1, \text{inl}(a_2))
	\end{array}$	
	\columnbreak
	$
	\begin{array}{rl}
	\equiv&\boldsymbol{\sigma}(\Delta(\eta(\phi_1,\phi_2)))(a_1, \text{inl}(a_2)) \\
	\equiv&\Delta(\eta(\phi_1,\phi_2))(a_1) (\text{inl}(a_2))\\
	\equiv&\boldsymbol{\kappa}_{a_1}(\phi_1(a_1),\phi_2(a_1))(\text{inl}(a_2))\\
	\equiv&\phi_1(a_1)(a_2)\;\;\,\equiv\text{inl}(a_1,a_2).
	\end{array}
	$
\end{multicols}
In the same way, $\Phi(\Psi(\text{inr}(b_1,b_2))):\equiv\text{inr}(b_1,b_2)$.
Then 
\[
\begin{calcu}
\expro{\prod\limits_{p:\sum_{x:A}P(x)+\sum_{x:A}Q(x)} \Phi(\Psi(p))=p}
\\
\explo{\simeq}{\!:\,$\boldsymbol{\kappa}$, ($\Pi$-range split rule)}
\\
\expro{\prod\limits_{a:\sum_{x:A}P(x)} \Phi(\Psi(\text{inl}(a)))=\text{inl}(a)
	\times \prod\limits_{b:\sum_{x:A}Q(x)} 
	\Phi(\Psi(\text{inr}(b)))=\text{inr}(b)}
\end{calcu}
\]
\[
\begin{calcu}
\explo{\simeq}{\!:\,$\boldsymbol{\sigma}_1\!\times\!\boldsymbol{\sigma}_2$,  {\sc eq}$_{\times}$ (\ref{EqProd})}
\\
\expro{\phantom{\times}\prod\limits_{a_1:A}\prod\limits_{a_2:P(a_1)} 
	\Phi(\Psi(\text{inl}(a_1,a_2)))=\text{inl}(a_1,a_2)\\
	&	\times \prod\limits_{b_1:A}\prod\limits_{a_2:Q(b_1)} 
	\Phi(\Psi(\text{inr}(b_1,b_2)))=\text{inr}(b_1,b_2)}
\\
\explo{\equiv}{Above computations}
\\
\expro{\phantom{\times}\prod\limits_{a_1:A}\prod\limits_{a_2:P(a_1)} 
	\text{inl}(a_1,a_2)=\text{inl}(a_1,a_2)\\
	&	\times \prod\limits_{b_1:A}\prod\limits_{a_2:Q(b_1)} 
	\text{inr}(b_1,b_2)=\text{inr}(b_1,b_2)}
\\
\explo{\stackrel{\mbox{\tiny $\wedge$}}{\mbox{\tiny :}}}{$u:\equiv\lambda a_1.\lambda a_2.\text{refl}_{\text{inl}(a_1,a_2)}$\,\,;\,\
	 $v:\equiv\lambda b_1.\lambda b_2.\text{refl}_{\text{inr}(b_1,b_2)}$}
\\
\expro{(u,v)}
\end{calcu}
\]
In the other direction,
\[
\begin{calcu}
\expro{\prod\limits_{p:\sum_{x:A}P(x)+Q(x)} \Psi(\Phi(p))=p}
\\
\explo{\simeq}{\!:\,$\boldsymbol{\sigma}$,\, $\Sigma$-consequent rule}
\\
\expro{\prod\limits_{x:A}\prod\limits_{y:P(x)+Q(x)} \Psi(\Phi(x,y))=(x,y)}
\\
\explo{\equiv}{Definition of $\Phi$}
\\
\expro{\prod\limits_{x:A}\prod\limits_{y:P(x)+Q(x)} 
	\Psi(\boldsymbol{\kappa}(\phi_1(x),\phi_2(x))(y))=(x,y)}
\\
\explo{\simeq}{\!:\,$\Delta$\,;\, $\varphi_x:\equiv \boldsymbol{\kappa}_x$\,;\,  $\Pi${\sc Eq}1}
\\
\expro{\prod\limits_{x:A}\phantom{\times} (\prod\limits_{w:P(x)} 
	\Psi(\boldsymbol{\kappa}(\phi_1(x),\phi_2(x))(\text{inl}(w)))=(x,\text{inl}
	(w)) \\
	& \phantom{\times\times}	\times \prod\limits_{z:Q(x)} 
	\Psi(\boldsymbol{\kappa}(\phi_1(x),\phi_2(x))(\text{inr}(z)))=(x,\text{inr}
	(z))) }
\\
\explo{\equiv}{Definition of $\boldsymbol{\kappa}$}
\\
\expro{\prod\limits_{x:A}\phantom{\times}( \prod\limits_{w:P(x)} 
	\Psi(\phi_1(x)(w))=(x,\text{inl}(w)) \\
	& \phantom{\times\times}	\times \prod\limits_{z:Q(x)} 
	\Psi(\phi_2(x)(z))=(x,\text{inr}(z))) }
\\
\explo{\equiv}{Definition of $\phi_1$ and $\phi_2$}
\\
\expro{\prod\limits_{x:A}\phantom{\times}( \prod\limits_{w:P(x)} 
	\Psi(\text{inl}(x,w))=(x,\text{inl}(w)) \\
	& \phantom{\times\times}	\times \prod\limits_{z:Q(x)} 
	\Psi(\text{inr}(x,z))=(x,\text{inr}(z))) }
\\
\explo{\equiv}{Definition of $\Psi$}
\\
\expro{\prod\limits_{x:A}\phantom{\times}( \prod\limits_{w:P(x)} 
	(x,\text{inl}(w))=(x,\text{inl}(w)) \\
	& \phantom{\times\times}	\times \prod\limits_{z:Q(x)} 
	(x,\text{inr}(z))=(x,\text{inr}(z))) }
\\
\explo{\stackrel{\mbox{\tiny $\wedge$}}{\mbox{\tiny :}}}{$u:\equiv \lambda x.(\lambda w.\text{refl}_{(x,\text{inl}(w))},\lambda z.\text{refl}_{(x,\text{inr}(z))})$}
\\
\expro{u}
\end{calcu}
\]
[\textbf{Translation}] rules correspond to the derived 
inference rules $\Pi${\sc Eq}2 and  $\Sigma${\sc Eq}2  which were proved in subsection \ref{LeibInf}\\[2mm]
[\textbf{Congruence}] rules correspond to the derived inference rules $\Pi${\sc Eq}1 and $\Sigma${\sc Eq}1 stated and proved in subsection 
\ref{LeibInf}\\[2mm]
[\textbf{Antecedent}] rules correspond to equivalences in first 
order logic that allow introducing the antecedent of an implication into the term of a logical operational when the quantified variables do not occur free in this antecedent. For HoTT, we only have an equivalence for the case of $\Pi$-types, $\Pi$-[\textbf{Antecedent}] rule. For $\Sigma$-types we have an equivalence only if the antecedent is a mere proposition. Namely,
\[
( P\rightarrow \prod\limits_{x:A}Q(x)) \; \simeq \; 
\prod\limits_{x:A} (P\rightarrow Q(x))<:
\]
and
\[\sum_{x:A}(P\rightarrow Q(x))
\rightarrow ( P\rightarrow \sum_{x:A}Q(x)) <:
\]
If $P\simeq \mathds 1<:$ then we get the equivalence.\\ [0.1cm]
The proof of $\Pi$-[\textbf{Antecedent}] rule appears in section \ref{InhArr}. We prove  $\Sigma$-[\textbf{Antecedent}] rule.
Let us consider the following deductive chain.
\[
\begin{calcu}
\expro{\sum\limits_{x:A}(P\rightarrow Q(x))\rightarrow (P\rightarrow 
	\sum\limits_{x:A}Q(x))}
\\
\explo{\simeq}{\!:\,$\boldsymbol{\sigma}$,\, $\Sigma$-consequent rule}
\\
\expro{\prod\limits_{x:A}((P\to Q(x))\rightarrow (P\to \sum\limits_{x:A}Q(x)))}
\\
\explo{\stackrel{\mbox{\tiny $\wedge$}}{\mbox{\tiny :}}}{$h(x)(u)(y):\equiv (x,u(y))$}
\\
\expro{h}
\end{calcu}
\]
This proves the first part. Now, If $P\simeq \mathds 1<:$, let 
\[
\psi:({\mathds 1}\to \sum_{x:A}Q(x))\rightarrow \sum_{x:A}({\mathds 1}\to Q(x))
\]
be defined by
\[
\psi(u):\equiv (\text{pr}_1(u(*)),\text{pr}_2\!\circ\! u).
\]

\section{Inhabiting arrows}\label{InhArr}
One of the tasks in homotopy type theory is to determine a formula for a 
function from type $A$ to a type $B$. We found that in several cases the 
structures of types $A$ and $B$ determine a natural matching of their objects 
defining  a function from $A$ to $B$. We call such a mapping a {\it 
canonical function}. An attempt to systematize this task is to precise the way 
in which we can get out of type $A$ through its eliminators and the way in 
which we can get in type $B$ through its constructors. To do so, we define the 
\textit{exit door} and the \textit{entry door} of a type. Of course, there will 
be types $A$ and $B$ for which there is no canonical function. This procedure is rather informal and has not relation with deductive chains, but allows us, in several cases, to find the canonical function.  

The entry door of a type is a $\lambda$-expression that represents a constructed 
object of the type, i.e., an  object of the type obtained from its 
constructors. The exit door of a type is a $\lambda$-expression that represents an 
eliminated object of the type, i.e., an object of the type constructed from the 
elimination of a generic object. For instance, the entry door of the type 
$\sum_{x:A}C(x)$ is the $\lambda$-expression
\[\left(u_1\!:\!A, u_2:C(u_1)\right) 
\]
because a constructed object of the type is a dependent pair of objects $u_1$ 
of 
type $A$ and $u_2$ of type $C(u_1)$. Then, we write
\[
\begin{calcu}
\expro{\sum\limits_{x:A}C(x)}
\\
\explo{\uparrow}{entry door}
\\
\expro{\left(\,u_1\!:\!A\,,\, u_2:C(u_1)\,\right)}
\end{calcu}
\]
The exit door of this type is the $\lambda$-expression 
\[\left(\, \text{pr}_1 (u)\!:\!A\,,\, \text{pr}_2(u)\!:\!C(\text{pr}_1(u)\,\right)
\]
because it is the dependent pair constructed from the elimination of a generic 
object $u$ of type $\sum\limits_{x:A}C(x)$ through their projections. We write
\[
\begin{calcu}
\expro{\left(\, \text{pr}_1 (u)\!:\!A\,,\, \text{pr}_2(u)\!:\!C(\text{pr}_1(u)\,\right)}
\\
\explo{\downarrow}{exit door}
\\
\expro{\sum\limits_{x:A}C(x).}
\end{calcu}
\]
The doors of a type can be used to determine a formula for a canonical function 
from a type to another, by matching the exit door of the source type  with the 
entry door of the destination type. For instance, let us determine a function 
from 
$\sum_{x:A}C(x)$ to itself.
This means that we have to determine an object $\Phi$ in the following link
\[
\begin{calcu}
\expro{\sum\limits_{x:A}C(x)}
\\
\explo{\leftarrow}{\!:\,$\Phi$}
\\
\expro{\sum\limits_{x:A}C(x),}
\end{calcu}
\]
i.e.  we have to match the exit door $\left(\, \text{pr}_1 (u)\!:\!A\,,\, 
\text{pr}_2(u)\!:\!C(\text{pr}_1(u)\,\right)$ and the entry door $\left(\Phi(u)_1\!:\!A, 
\Phi(u)_2:C(\Phi(u)_1)\right)$ of the type $\sum_{x:A}C(x)$, task that 
we represent with the following matching diagram  
\[
\begin{calcu}
\expro{\sum\limits_{x:A}C(x)}
\\
\explo{\uparrow}{entry door}
\\
\expro{\left(\Phi(u)_1\!:\!A, \Phi(u)_2:C(\Phi(u)_1)\right)}
\\
\explo{\mapsfrom}{Looking for definition}
\\
\expro{\left(\, \text{pr}_1 (u)\!:\!A\,,\, \text{pr}_2(u)\!:\!C(\text{pr}_1(u))\,\right)}
\\
\explo{\downarrow}{exit door}
\\
\expro{\sum\limits_{x:A}C(x),}
\end{calcu}
\]
where $\mapsfrom$ means that some sort of symbolic matching between two 
expressions must be discovered. By matching the doors we get 
\[
\Phi(u):\equiv (\text{pr}_1(u), \text{pr}_2(u)).
\]
Observe that the canonical function in this case is not the identity 
function.\\ 
[0.1cm]
Let us determine the canonical function $\Phi$ from $\prod_{x:A}B(x)$ to 
itself. The corresponding matching diagram is
\[
\begin{calcu}
\expro{\prod\limits_{x:A}B(x)}
\\
\explo{\uparrow}{entry door}
\\
\expro{\lambda(x\!:\!A).(\Phi(f)\!:\!B)}
\\
\explo{\mapsfrom}{?}
\\
\expro{\lambda(x\!:\!A).(f(x)\!:\!B(x))}
\\
\explo{\downarrow}{exit door}
\\
\expro{\prod\limits_{x:A}B(x).}
\end{calcu}
\]
Therefore, by matching, we get
\[
\Phi(f)(x):\equiv f(x).
\]
which, by uniqueness, is the identity function. \\ [0.1cm]
We now present some examples illustrating this technique.\\[0.1cm]
\textbf{$\Pi$-distribution over arrows}. As promised in section \ref{DistArr}, we show how to obtain the canonical function $\Phi:\equiv \lambda u.\Phi(u)$ of the type
\[\prod\limits_{x\!:\!A} (P(x) \rightarrow Q(x)) 
\rightarrow
( \prod\limits_{x\!:\!A} P (x) \rightarrow \prod\limits_{x: A} Q (x)) . 
\]
For that, the corresponding entrance and exit doors are made to coincide
\[
\begin{calcu}
\expro{\prod\limits_{x: A} P (x) \rightarrow \prod\limits_{x: A} Q (x)}
\\
\explo{\uparrow}{entry door}
\\
\expro{\lambda (z\!:\!\prod\limits_{x: A} P (x)). \lambda (x\!:\!A). 
\Phi(u)(z)(x)}
\\
\explo{\mapsfrom}{?}
\\
\expro{\lambda (x\!:\!A).\lambda (y\!:\!P(x)).u(x)(y)}
\\
\explo{\downarrow}{exit door}
\\
\expro{\prod\limits _{x:A}(P(x) \rightarrow Q(x))}
\end{calcu}
\]
obtaining
\[\Phi(u)(z)(x):\equiv u(x)(z(x)).
\]
$\Pi$-[\textbf{Antecedent}] rule. In order to prove that
\[
( P\rightarrow \prod\limits_{x:A}Q(x) ) \; \simeq \;
\prod\limits_{x:A} (P\rightarrow Q(x))<:
\]
we have to determine a 4-tuple $(\Phi, \Phi',\alpha, \alpha')$ inhabiting the equivalence type.  
Consider the following entry-exit door arguments:
\[
\begin{calcu}
\expro{P\rightarrow \prod\limits_{x:A}Q(x)}
\\
\explo{\uparrow}{entry door}
\\
\expro{\lambda(y\!:\!P).\lambda(x\!:\!A).(\Phi(u)(y)(x):Q(x) )}
\\
\explo{\mapsfrom}{$\Phi(u)(y)(x):\equiv u(x)(y)$}
\\
\expro{\lambda(x\!:\!A).\lambda(y\!:\!P).(u(x)(y):Q(x)) }
\\
\explo{\downarrow}{exit door}
\\
\expro{\prod\limits_{x:A} (P\rightarrow Q(x)),}
\end{calcu}
\]
and
\[
\begin{calcu}
\expro{\prod\limits_{x:A} (P\rightarrow Q(x))}
\\
\explo{\uparrow}{entry door}
\\
\expro{\lambda(x\!:\!A).\lambda(y\!:\!P).(\Phi'(v)(x)(y)\!:\!Q(x)) }
\\
\explo{\mapsfrom}{$\Phi'(v)(x)(y):\equiv v(y)(x)$}
\\
\expro{\lambda(y\!:\!P).\lambda(x\!:\!A).(v(y)(x)\!:\!Q(x))}
\\
\explo{\downarrow}{exit door}
\\
\expro{P\rightarrow \prod\limits_{x:A}Q(x).}
\end{calcu}
\]
Observe that, by definition of $\Phi$ and $\Phi'$, 
\[\Phi'(\Phi(u))(x)(y)\equiv \Phi(u)(y)(x)\equiv u(x)(y)\]
and
\[\Phi(\Phi'(v))(y)(x)\equiv\Phi'(v)(x)(y)\equiv v(y)(x).\]
This shows that $\Phi'$ and $\Phi$ are each other inverses, and then, that   $\Phi'\circ\Phi\sim \text{id}<:$ and $\Phi\circ\Phi'\sim \text{id}<:$

\section{Conclusions}
We were able to obtain a formal deduction method in HoTT based on deduction chains; and found that the most important equational axioms and rules of a calculation version of intuitionistic logic (ICL) have a counterpart as derivable judgments in HoTT. Some of this judgments correspond to homotopic equivalence versions of the induction operators of basic types in HoTT.

We think that the use of deductive chains to formally prove HoTT theorems, in comparison with rigorous proofs written on paper by a human, is more effective, clear and readable. This is so, because the proofs are made of formally precise linearly chained modules which characterize the linear proof formats we call deductive chains. This way of proving, in our view, has the advantage of, on one hand, preserve formality avoiding ambiguities and imprecisions that may come with rigorous but colloquial proofs typical of the working mathematician; and on the other hand, they are constructed via very simple and precise steps, amenable to be made by hand. We hope to have helped demythify the wide belief that formal proofs are messy and very long to be readable and performable, in a practical way, by humans.\\
This work, appears to make possible the restatement of the whole HoTT in terms of an appropriate calculus of equational deduction.\\ 
Finally, we expect that our research will motivate exploring the proof theory associated to calculational methods of proof. We also think that it would be worthwhile to develop proof assistants and verifiers to support the automation of these methods.
\bibliographystyle{abbrv}
\bibliography{DeductiveChains}

\begin{thebibliography}{10}

\bibitem{AAB17}
E.~Acosta, B.~Aldana, J.~Boh{\'o}rquez, and C.~Rocha.
\newblock Axiomatic set theory {\`a} la {D}ijkstra and {S}cholten.
\newblock In A.~Solano and H.~Ordo{\~{n}}ez, editors, {\em Advances in
  Computing}, pages 775--791, Cham, 2017. Springer International Publishing.

\bibitem{Bac03}
R.~Backhouse.
\newblock {\em Program Construction: Calculating Implementations from
  Specifications}.
\newblock John Wiley and Sons, Inc., 2003.

\bibitem{Bar02a}
H.~Barendregt and E.~Barendsen.
\newblock Autarkic computations in formal proofs.
\newblock {\em J. Automated Reasoning}, 28(3):321--336, 2002.

\bibitem{Boh05}
J.~Boh{\'o}rquez and C.~Rocha.
\newblock Towards the effective use of formal logic in the teaching of discrete
  math.
\newblock 6th International Conference on Information Technology Based Higher
  Education and Training. ITHET., 2005.

\bibitem{Boh08}
J.~A. Boh{\'o}rquez.
\newblock Intuitionistic logic according to {D}ijkstra's calculus of equational
  deduction.
\newblock {\em Notre Dame J. Form. Log.}, 49(4):361--384, 2008.

\bibitem{Boh15}
J.~A. Bohorquez.
\newblock Calculational solutions to combinatorial problems.
\newblock 10th Computing Colombian Conference (10CCC), 2015.

\bibitem{DS90}
E.~W. Dijkstra and C.~S. Scholten.
\newblock {\em Predicate Calculus and Program Semantics}.
\newblock Springer Verlag, 1990.

\bibitem{FvG99}
W.~H.~J. Feijen and A.~J.~M. van Gasteren.
\newblock {\em On a method of multiprogramming}.
\newblock Springer-Verlag New York, Inc., New York, NY, USA, 1999.

\bibitem{DBLP:conf/procomet/gries98}
D.~Gries.
\newblock Teaching calculational logic.
\newblock In D.~Gries and W.~P. de~Roever, editors, {\em PROCOMET}, volume 125
  of {\em IFIP Conference Proceedings}, pages 9--10. Chapman {\&} Hall, 1998.

\bibitem{Gri93}
D.~Gries and F.~B. Schneider.
\newblock {\em A Logical Approach to Discrete Math}.
\newblock Texts and Monographs in Computer Science. Springer Verlag, 1993.

\bibitem{Lif01}
V.~Lifschitz.
\newblock On calculational proofs.
\newblock {\em Ann. Pure Appl. Logic}, 113(1-3):207--224, 2001.

\bibitem{Misra01}
J.~Misra.
\newblock {\em A Discipline of Multiprogramming: Programming Theory for
  Distributed Aplications}.
\newblock Monographs in Computer Science. Springer-Verlag, New York, 2001.

\bibitem{hottbook}
T.~{Univalent Foundations Program}.
\newblock {\em Homotopy Type Theory: Univalent Foundations of Mathematics URL
  https://homotopytypetheory.org/book}.
\newblock Institute for Advanced Study, 2013.

\bibitem{vGast90}
A.~J.~M. van Gasteren.
\newblock {\em On the Shape of Mathematical Arguments}, volume 445 of {\em
  Lecture Notes in Computer Science}.
\newblock Springer, 1990.

\end{thebibliography}

\end{document}